\documentclass[preprint,12pt]{elsarticle}

\usepackage{amssymb}
\usepackage{amsthm}
\usepackage{amsmath}
\usepackage{color}
\usepackage[labelfont=bf,list=true,font=footnotesize]{subcaption}
\newdefinition{rmk}{Remark}

\usepackage{lineno}

\journal{Elsevier}

\begin{document}
\begin{frontmatter}



\title{A multi-model study of the air pollution related to traffic flow in a two-dimensional porous metropolitan area}


\author[guadalajara]{N. Garc\'\i a-Chan}
\ead{nestorg.chan@cucei.udg.mx}
\address[guadalajara]{Depto. F\'\i sica. C.U.~Ciencias Exactas e Ingenier\'\i as.
Universidad de Guadalajara.\\ 44420 Guadalajara. Mexico}
\author[vigo]{L.J. Alvarez-V\'azquez\corref{cor1}{}}
\ead{lino@dma.uvigo.es}
\address[vigo]{Depto.~Matem\'atica Aplicada II. E.I.~Telecomunicaci\'on.
Universidad de Vigo. \\ 36310 Vigo. Spain}
\author[vigo]{A. Mart\'\i nez}
\ead{aurea@dma.uvigo.es}
\author[lugo]{M.E.~V\'azquez-M\'endez}
\ead{miguelernesto.vazquez@usc.es}
\address[lugo]{Depto. Matem\'atica Aplicada. E.P.S. Universidad de
Santiago de Compostela. \\ 27002 Lugo. Spain}         

\begin{abstract}
In this paper, a useful reinterpretation of the city as a porous medium justifies the application of well-known models on fluid dynamics to develop a multi-model study of urban air pollution due to traffic flow in a large city. Thus, to simulate the traffic flow through the city we use a nonconservative macroscopic traffic model combining the continuity equation with the Darcy-Brinkman-Forchheimer equations.  For the air flow,  regarding the emission rate of
$\rm{CO_2}$ and its dispersion in the atmosphere, we combine a microscopic model -based on regression techniques but depending on vehicles' velocity and acceleration- with a classical convection-diffusion-reaction transport model.  To solve numerically above PDEs models, the finite element method of Lagrange $\rm{P_1}$ type along with suitable time marching schemes (like the strong stability preserving scheme) were sufficient to obtain stable numerical solutions. Several computational tests were run on a realistic scenario inspired by the Metropolitan Area of Guadalajara (Mexico),  showing not only the influence of the urban landscape (that is,  the porosity) on traffic flow,  air flow,  and pollution transport,  but also other interesting phenomena such as rarefaction traffic waves.
\end{abstract}




\end{frontmatter}


\section{Introduction}
In large urban areas, traffic flow is the primary source of air pollution; however, studying this phenomenon is challenging due to various factors. The high number of vehicles (thousands or millions) contributes to the complexity of the problem, as do traffic jams resulting from poor urban planning, congestion and accidents.  The intricate urban landscape, characterized by a complex network of primary and secondary avenues and streets,  the wide range of available fuels (gas, diesel), and the diverse types of vehicles (combustion cars and trucks, hybrid and electric vehicles, heavy-duty vehicles) further complicate the situation.  However, despite the apparent chaos, the traffic flow exhibits well-observed patterns, including attraction zones (typically in downtown areas and commercial zones), periods of peak and off-peak hours, and drivers' preferences for primary roads. The integration of these elements needs the development of sophisticated mathematical models, particularly the characterization of the urban landscape in a metropolis of considerable size. In this context, the authors recently proposed a novel model that interprets the large city as a porous medium, with porosity serving as a metric for the proportion of available space for streets and the complementary space occupied by buildings \cite{garcia-chan2025}. This approach enables the formulation of a nonconservative traffic flow model consisting of a continuity equation based on a source/sink dynamic (traffic demand/off-street parking), a momentum equation of Darcy-Brinkman-Forchheimer type with a relaxation term, and a pressure correction. This model was employed with success to obtain cars density and local speed.   
 
A further critical element that must be incorporated into the model is the emissions rate released by traffic flow.  The extensive scientific literature on the topic delineates two alternative standpoints \cite{Jiang2016}: the macroscopic models and the microscopic ones.  The first approach involves the interpolation of mean values of various pollutants measured on specific streets and avenues of the road network, extending these measurements to the entire network while neglecting key microscopic factors such as vehicle velocity and acceleration \cite{CARB2007,EPA2003,CopertIII}.  However, in a two-dimensional study, this kind of model encompasses the entire city losing spatial resolution.  In contrast, microscopic emission rate models consider vehicle velocity and acceleration, using a regression polynomial whose coefficients are adjusted based on the type of combustible and the pollutant of interest.  Consequently, it could be profitable to convert microscopic emission rates into macroscopic ones when the distribution of vehicle density, velocity, and acceleration is available at any city location (see \cite{Jiang2016}).

The spatial availability for air circulation influences the air flow dynamics, thus, in a wide and large street surrounded by buildings the well-known phenomenon of the wind cannon is produced. This phenomenon is characterized by sharp high wind velocities along the street. In contrast, short and narrow streets generate air stagnation with almost zero velocity. It is important to observe that the space available for air circulation is complement by the space occupied by buildings, this drives us -again- to a porous medium with its fluid and solid parts. 
This air behavior suggests the use of a porous media model for air flow as well as for the traffic flow (see \cite{hu2012} for a stationary case). The authors already published a study on the flow of air as an incompressible gas through an urban porous city, so we take this model in our study \cite{garcia-chan2023,garcia-chan2024}.
Concerning the distribution of the pollutant in the atmospheric boundary layer, wind flow and diffusion effects transport the pollutant until reaching remote zones from the release point. Also, depending on the pollutant type, its concentrations could decay by kinetic reactions or by gravitational settling. The pollutant transport is then modeled by a classical advection-diffusion-reaction PDE (see \cite{sportisse2010} and the references therein).  

In this work, a multi-model consisting of models based on fluid dynamics in porous media combined with a microscopic emissions model and a classical dispersion-diffusion-reaction model is used to address the main elements of the urban air pollution phenomenon, resulting in a comprehensive numerical study.  The paper is organized as follows. In Section 2, the three mathematical models are formulated, and their most important elements are explained.  In Section 3, we briefly explain how to solve the models combining the $\rm{P_1}$ Lagrange-type finite element method with suitable marching time schemes, along with the references of previous authors' papers where the reader can find further details.  In Section 4, the results of several numerical experiments are shown, evaluating the influence of the urban landscape on traffic and air flow, which influences the $\rm{CO_2}$ emission rate and its dispersion. Finally, in Section 5, we state several claims supported by our numerical experiences regarding air pollution due to traffic flow in large urban areas.

%
\section{The mathematical models}\label{sec:mathmodel}

In this section, we proceed to formulate the mathematical models involved in our meta-model study, wherein fluid dynamics on porous media plays a pivotal role along with a regression polynomial model to estimate the emissions rate of carbon dioxide. It is pertinent to note that, with the exception of the microscopic emissions rate model, all models have been analyzed by the authors in their recent publisher papers. Consequently, only its most relevant elements are exhibited here, and interested readers are encouraged to consult the references for further details.

A remarkable strength of this paper relies in the fact that that the same domain (i.e., the same mesh) may be used for all the models despite their differences in the definition of the corresponding domains. Thus, we denote our domain as $\Omega$, which contains the rural zone $\Omega_{r}$ and the urban zone $\Omega_{u}$,  surrounding a few inner, natural obstacles that are not included into the domain $\Omega_u$.  Domain boundary is denoted by $\Gamma$, and is divided into the following disjoint line segments: the obstacles' walls $\Gamma_{w}$,  the inlet airflow $\Gamma_{in}$,  and the outlet airflow $\Gamma_{out}$.  A particular example of domain $\Omega$ can be seen in below Figure (\ref{fig1:GDL-mesh}) in Section 4.  So let us start with the model for the cars density and its velocity. 

\subsection{A nonconservative macroscopic traffic flow model in an urban-porous medium}

We use the traffic flow model formulated in \cite{garcia-chan2025} which considers some key elements such as viscosity, Darcy law and its corrections, and uses a desired direction to conduct to drivers in the momentum equation, meanwhile, the continuity equation is formulated considering a nonconservative dynamic which includes the diffusion of cars, and a source/sink dynamic given by traffic-demand and street-off parking terms. The desired speed is computed from an Eikonal problem,  whose formulation was improved in \cite{garcia-chan2025},  based on the one formulated in \cite{Jiang2016} (see equation (\ref{mod:eikonal})).  

Then, we seek for the cars density $\rho\ [{\rm veh/km^2}]$ and their local traffic speed $\mathbf{u}=(u_1,u_2)\ [{\rm km/h}]$ such that satisfy the following system of PDEs, where the dot ( $\dot{ }$ ) denotes the time derivative and vector $\nabla$ stands for the classical spatial gradient operator: 
\begin{subequations}\label{mod:traffic-flow}
\begin{align}
   \epsilon\dot{\rho} + \nabla\cdot (\rho\mathbf{u}) - \nabla\cdot (\epsilon\nu\nabla\rho) + \epsilon\kappa\rho = (1-\epsilon)q &\quad \mbox{ in }\Omega\times (0,T),\label{eq:trafficflow-continuity}\\
   \dot{\mathbf{u}} + \frac{1}{\epsilon} (\mathbf{u}\cdot\nabla) \mathbf{u} - c^2\rho(\nabla\cdot\mathbf{u})\mathbf{1} - \frac{\mathbf{u}_d - \mathbf{u}}{\tau} 
   & \nonumber\\
 = \frac{\epsilon}{\rho}(\nabla\cdot\frac{\mu}{\epsilon}\nabla) \mathbf{u}  - \frac{\epsilon \mu}{\rho K} \mathbf{u} - \frac{\epsilon F}{\sqrt{K}}\Vert\mathbf{u}\Vert \mathbf{u} & \quad \mbox{ in }\Omega\times (0,T),\label{eq:traffic-flow-momentum}   
\end{align}
\end{subequations}
where $\epsilon$ is the urban porosity, $(1-\epsilon)q\ [{\rm veh/km}^2{\rm /h}]$ is the traffic demand function associated to the building-blocks (solid-phase), $\mathbf{u}_d\ [{\rm km/h}]$ is the desired speed, $c^2$ is the traffic sonic speed (computed as $c^2 = \theta/\rho_{max}$ with $\theta$ an anticipation coefficient and $\rho_{max}$ a fixed road capacity), $K$ is the urban permeability, $F$ is the Forchheimer coefficient, $\nu\ [{\rm km}^2{\rm /h]}$ is the coefficient of mass diffusion, $ \kappa\ [{\rm h}^{-1}]$ is the absorption (parking) rate from the streets (fluid-phase) to parking-spaces inside building-blocks (solid-phase), $\tau\ [{\rm h}]$ is the relaxation time, and $\mu\ [{\rm km}^2{\rm /h}]$ is the viscosity. 

The system (\ref{mod:traffic-flow}) is complemented with suitable initial and boundary conditions. So, we impose a Neumann homogeneous condition, and a slip-type condition on the two boundary segments to avoid escaping cars by diffusion, and to enforce a tangential flow of vehicles when an obstacle is on their path:
\begin{subequations}\label{bnd:traffic-flow}
\begin{align}
\nabla\rho\cdot\mathbf{n} = 0 &\quad \mbox{ on } (\Gamma_{w}\cup\Gamma_{in}\cup\Gamma_{out}) \times (0,T), \label{bnd:traffic-flow-wall-bnd}\\   
\mathbf{u}\cdot\mathbf{n} = 0 &\quad \mbox{ on } \Gamma_{w} \times (0,T),  \label{bnd:traffic-flow-slip}\\
\nabla\mathbf{u} \mathbf{n} = \mathbf{0} &\quad \mbox{ on } (\Gamma_{in}\cup\Gamma_{out}) \times (0,T).\label{bnd:traffic-flow-neumann-bnd}
\end{align}
\end{subequations}  
(As well known,  last condition $\nabla\mathbf{u} \mathbf{n} = \mathbf{0}$ is equivalent to separate conditions $\nabla u_1\cdot\mathbf{n} = \nabla u_2\cdot\mathbf{n} = 0$. )

We also assume that, at the initial time $t=0$, traffic density presents a given distribution $\rho^0$ where all cars are at rest, that is:
\begin{subequations}\label{cond-ini:traffic-flow}
\begin{align}
\rho(.,0) = \rho^0(.) & \quad \mbox{ in } \Omega \label{cond-ini:traffic-flow-density},\\
\mathbf{u}(.,0) = \mathbf{0} & \quad \mbox{ in } \Omega\label{cond-ini:traffic-flow-velocity}.
\end{align}
\end{subequations}
\subsection{Characterization of the urban and rural zones, the desired direction, and the travel cost}

We must recall here that all the parameters in our model are continuous in the domain, in spite of the domain also presents clearly defined urban and rural zones.  According to this, we must define the local travel cost and the desired traffic speed generated from it also as continuous functions.   So, let $1_{\Omega_{u}}$ and $1_{\Omega_{r}}$ be the indicator functions for the urban zone $\Omega_u$ and the rural zone $\Omega_r$,  respectively.  That is,  $1_{\Omega_{u}}(\mathbf{x}) = 1$ ($1_{\Omega_{r}}(\mathbf{x}) = 0$) if $\mathbf{x}\in \Omega_{u}$,  and $1_{\Omega_{u}}(\mathbf{x}) = 0$ ($1_{\Omega_{r}}(\mathbf{x}) = 1$) if $\mathbf{x}\in\Omega_r$. Therefore,  using the above indicator functions we formulate the expression for the local travel cost as: 
\begin{equation}
f(\rho) = U_{max}(1 - 1_{\Omega_{u}}\frac{\rho}{\rho_{max}} - 1_{\Omega_{r}}) + 1_{\Omega_{r}}.
\end{equation}
In this way, if drivers flow inside the city,  the cost $f$ is reduced to the classical relation $U_{max}(1 - \rho/\rho_{max})$ where its value depends on the cars density \cite{garcia-chan2025},  being its maximum value $f=U_{max}$ when $\rho=0$,  and its minimum value $f=0$ when $\rho = \rho_{max}$.  In the other way,  if drivers go from the city to the rural zone, the travel cost is the unity ($f=1$),  which does not depend on the density.  Thus,  since the desired speed is computed using the inverse $1/f$ (see equation (\ref{eq:eikonal})),  driving inside the city has a cheaper cost than driving in rural roads.  Now,  once defined the travel cost, we can use an Eikonal approach to find the desired speed, which is written here as in \cite{Jiang2016}:
\begin{subequations}\label{mod:eikonal}
\begin{align}
\Vert\varphi\Vert = 1/f(\rho) & \quad \mbox{ in } \Omega, \label{eq:eikonal}\\
\varphi = 0 & \quad \mbox{ on }\Gamma_{exit}, \label{bnd:eikonal-exit}\\
\nabla\varphi\cdot\mathbf{n} & \quad \mbox{ on }\Gamma_{w}\cup\Gamma_{in}\cup\Gamma_{out},  \label{bnd:eikonal-wall-boundary}
\end{align}
\end{subequations}
being the potential $\varphi$ the instantaneous travel cost depending on the inverse of the transportation cost.  It is noteworthy that problem (\ref{mod:eikonal}) models a direction $\nabla\varphi$ pointing towards an ``exit'' where $\varphi=0$ (located on the boundary segment $\Gamma_{exit}$), but circumventing inner obstacles and avoiding crossing the domain limits.  Also, model (\ref{mod:eikonal}) is nonlinear, but can be turned into a linear problem following the steps given in \cite{axthelm2016,churbanov2019}.  
Specifically,  for our case,  we change the concept of ``exit'' by an attraction point enforced by a function $G$, use the linearization based on the change of variable $\psi = e^{-\varphi / \eta}$, and adapt our conditions on the boundary segments $\Gamma_{in}, \Gamma_{out}$ and $\Gamma_w$. After all these steps we can reformulate the problem (\ref{mod:eikonal}) as (see \cite{garcia-chan2025} for further details):
\begin{subequations}\label{mod:eikonallineal}
\begin{align}
\eta^2\Delta\psi - \psi/f^2(\rho) = G & \quad \mbox{ in }\Omega,\label{eq:eikonal-lineal}\\
\nabla\psi\cdot\mathbf{n} = 0 & \quad \mbox{ on } \Gamma_w\cup\Gamma_{in}\cup\Gamma_{out},\label{bnd:eikonal-lineal-neumann}
\end{align}
\end{subequations}
with $\eta>0$ a regularization parameter and $G$ the function related to drivers' destination.  Obviously,  the solution $\psi$ of the linear problem (\ref{mod:eikonallineal}) is related to the potential $\varphi$ by the relation $\varphi = -\eta \ln(\psi)$.  Thus, we finally define the desired speed as:
\begin{equation}
\mathbf{u}_d(\rho) = -f(\rho)\nabla\varphi/\Vert\nabla\varphi\Vert.
\end{equation}

\subsection{Evaluating $\rm{CO_2}$ emission rates using a microscopic emission model}

The computation of emission rates of carbon dioxide $\rm{CO_2}$ in this paper is achieved through the utilization of a classical microscopic emission model, which is formulated as a regression polynomial,  as can be found in \cite{panis}.  This model needs the use of scalar functions for microscopic velocity and acceleration (in $[\rm{m/s}]$ and $[\rm{m/s^2}]$,  respectively).  In the context of this model, the emission rate of $\rm{CO_2}$ (in $[\rm{g/s/veh}]$) is derived as a function of these microscopic velocity $U$ and acceleration $a$: 
\begin{equation}\label{eq:instant-emissions}
E_{\rm CO_2} = \max\left\{f_1 + f_2 U + f_3 U^2 + f_4 a + f_5 a^2 + f_6 Ua, 0.0\right\}  
\end{equation}
where coefficients $f_i, \ i=1,2,\ldots,6,$ have suitable values for $\rm{CO_{2}}$.  Many important quantities on air quality can computed from equation (\ref{eq:instant-emissions}).  In our case,  the emission concentration $EC(\mathbf{x},t)$ (in $[\rm{kg/km^2/h}]$) will be given by:
\begin{equation}\label{eq:concent-emissions}
EC = \rho\,E_{\rm CO_2}\gamma_1
\end{equation}
being $\rho$ the cars density,  and $\gamma_1=3.6$ the factor to adjust microscopic dimensions $[\rm{g/s}]$ into macroscopic ones $[\rm{kg/h}]$.

Therefore, we need compute $U$ and $a$ from our macroscopic velocity $\mathbf{u}$ and acceleration $\mathbf{a}$. The first one is straightforward computed with the Euclidean norm of the traffic flow velocity, adjusted by a conversion parameter from $[\rm{km/h}]$ to $[\rm{m/s}]$ \cite{Jiang2016}:
\begin{equation}\label{eq:scalar-velocity}
U = \Vert\mathbf{u}\Vert/\gamma_1.
\end{equation}   
For the latter,  we use the definition of the acceleration as the material derivative $\mathbf{a} = D\mathbf{u}/Dt = \dot{\mathbf{u}} + \frac{1}{\epsilon} (\mathbf{u}\cdot\nabla) \mathbf{u}$ of fluid dynamics in porous media (see \cite{das2018_book}) from a Lagrangian perspective.  In this regard, given the traffic flow velocity $\mathbf{u}$, the desired speed $\mathbf{u}_d$ and the cars density $\rho$,  the vectorial acceleration can be computed from equation (\ref{eq:traffic-flow-momentum}) with the following formula:  
\begin{equation}\label{eq:acceleration-vector}
\mathbf{a} = \frac{\mathbf{u}_d - \mathbf{u}}{\tau} + c^2\rho(\nabla\cdot\mathbf{u})\mathbf{1} + \frac{\epsilon}{\rho}(\nabla\cdot\frac{\mu}{\epsilon}\nabla) \mathbf{u}  - \frac{\epsilon \mu}{\rho K} \mathbf{u} - \frac{\epsilon F}{\sqrt{K}}\Vert\mathbf{u}\Vert \mathbf{u}. 
\end{equation}
Nevertheless, a second step must be done to pass from the vector field $\mathbf{a}$ to a scalar function.  Thus, in agreement with \cite{Jiang2016},  we use the following equation:  
\begin{equation}\label{eq:scalar-acceleration}
a = (\mathbf{a}\cdot \mathbf{u}/U)\gamma_2,
\end{equation}
where $\gamma_2=10^{-3}/(3.6)^2$ is a parameter to escalate the macroscopic units $[\rm{km/h^2}]$ into microscopic dimensions $[\rm{m/s^2}]$. 
\subsection{The air flow model}
The urban landscape with its buildings and streets influences the available free space for air circulation.  So, wide and large streets with large buildings on their sides generate an urban cannon phenomenon characterized by high wind velocities.  In contrast, narrow streets with a variety of building sizes can act as barriers to wind circulation. Therefore, as the traffic flow, computational fluid dynamics models on porous media are proposed to model the air flow inside the city (see \cite{garcia-chan2023,garcia-chan2024}).  In order to find the average air velocity $\mathbf{v}$ it is necessary to solve the following system of PDEs:
\begin{subequations}\label{mod:airflow}
\begin{align}
\dot{\mathbf{v}} + \frac{1}{\epsilon}(\mathbf{v\cdot\nabla})\mathbf{v} = -\frac{\epsilon}{\rho_a}\nabla P - \frac{\epsilon\mu}{\rho_a K}\mathbf{v} + \frac{\epsilon}{\rho_a}(\nabla\cdot\frac{\mu}{\epsilon}\nabla)\mathbf{v} 
\nonumber\\ - \frac{\epsilon C_F}{\sqrt{K}}\Vert\mathbf{v}\Vert\mathbf{v} & \quad \mbox{ in } \Omega\times(0,T), \label{eq:airflow-momentum}\\ 
\nabla\cdot\mathbf{v} = 0 & \quad \mbox{ in } \Omega\times(0,T), \label{eq:airflow-continuity} 
\end{align}
\end{subequations}
where $\rho_a$ [$\rm{kg/km^3}$] is the air density, and $\mu_a$ [$\rm{km^2/h}$] is the dynamical viscosity of air. Due to the fact that the urban landscape is the same for the traffic flow and the air flow, the porosity $\epsilon$,  the permeability $K$, and the Forchheimer coefficient $C_F$ remain the same in both models.  We also assume that air is an incompressible fluid due to the fact that,  since energy interchange (air temperature) is not included in our study,  air density will remain constant. 
To define the initial/boundary conditions to complement the model (\ref{mod:airflow}), we assume the following statements: air is calm at the initial time,  the air flow through the inlet boundary segment is given by $\mathbf{v}_{in}$, the wind slips around obstacles, there are no changes on wind direction by the wind itself and the pressure on the outlet boundary segment, and pressure does not influence on walls and inlet boundary segments.  In agreement with all these assumptions we impose the following conditions:
\begin{subequations}\label{bnd:airflow}
\begin{align}
\mathbf{v} = \mathbf{v}_{in} & \quad  \mbox{ on }\Gamma_{in} \times (0,T), \label{bnd:airflow-inlet}\\
\mathbf{v}\cdot\mathbf{n} = 0 & \quad   \mbox{ on }\Gamma_{w}  \times (0,T), \label{bnd:airflow-slip-wall}\\
\frac{\mu_a}{\epsilon}\nabla\mathbf{v} \mathbf{n} - P \mathbf{n} = \mathbf{0} & \quad  \mbox{ on } \Gamma_{out} \times (0,T), \label{bnd:airflow-outlet}\\
\nabla P \cdot\mathbf{n}= 0 & \quad  \mbox{ on } (\Gamma_{in}\cup\Gamma_{w}) \times (0,T),  \label{bnd:airflow-pressure} \\
\mathbf{v}(.,0) =  \mathbf{0} & \quad \mbox{ in }\Omega.\label{ini:airflow}
\end{align}
\end{subequations}

\begin{rmk} We must distinguish the average velocity $\mathbf{v}$ on a reference elementary volume (REV) from the local velocity $\tilde{\mathbf{v}}$ on the REV (which is the velocity computed, for instance,  in \cite{garcia-chan2023,garcia-chan2024}).  Both velocities are related by the equality \cite{das2018_book}:
\begin{equation}
\mathbf{v} = \frac{\tilde{\mathbf{v}}}{\epsilon}.
\end{equation} 
However, in this paper we will compute directly the average velocity $\mathbf{v}$ due to the numerical schemes chosen here,  that allow its straightforward computation.  
\end{rmk}
\subsection{The air pollution model}
Macroscopic traffic flow at the entire city generates instantaneous emissions at a microscopic level which,  after translated into macroscopic emission concentration $EC$,  is then transported by wind and diffusion through the city,  evacuating air pollution in favor of inhabitants' health.  Also, pollutants can suffer deposition effects by gravity action or by low order reactions. This leads to the formulation of a classical advection-diffusion-reaction PDE to model air-pollution transport \cite{sportisse2010},  but including porosity since transport is done -again- through the streets and avenues (fluid-phase) of the porous city.  Thus, let $\phi$ be the air pollutant concentration of $\rm CO_2$ (in $[\rm{kg/km^2}]$) such that satisfies the air pollution transport model:
\begin{subequations}\label{mod:airpollution}
\begin{align}
\epsilon\dot{\phi} + \mathbf{v}\cdot\nabla\phi - \nabla\cdot\epsilon\mu_{\phi}\nabla\phi + \epsilon\sigma\phi = EC & \quad \mbox{ in }\Omega\times(0,T) , \label{eq:pollution}\\
\epsilon\mu_{\phi}\nabla\phi \cdot \mathbf{n} - \phi \mathbf{v}\cdot\mathbf{n} = 0 & \quad \mbox{ on }\Gamma_{in}\times(0,T), \label{bnd:pollution-inlet}\\
\epsilon\mu_{\phi}\nabla\phi \cdot \mathbf{n} = 0 &\quad  \mbox{ on }(\Gamma_{out}\cup\Gamma_{w})\times(0,T),  \label{bnd:pollution-out-wall}\\
\phi(.,0) = \phi^0(.) & \quad \mbox{ in } \Omega, \label{inicond:pollution} 
\end{align}
\end{subequations}
where $\mathbf{v}$ is the air flow velocity from model (\ref{mod:airflow}), $\mu_{\phi}$ is the diffusivity coefficient, $\sigma>0$ denotes the deposition rate,  and $EC$ is the macroscopic emission concentration of $\rm CO_2$ computed from equation (\ref{eq:concent-emissions}).  Above formulated initial/ boundary conditions, indicate that air pollution over domain $\Omega$ is initially given by a known distribution $\phi^0$,  that the flow of $\rm CO_2$ is null through the boundary inlet segment, and that pollutant diffusion on obstacle walls and outlet boundary is forbidden.  Also, it must be noted the influence of porosity (i.e., building density) on the time variability, the diffusion term, and the reaction term.

Due to the large size of metropolitan areas, its irregular distribution of buildings and streets, the dynamics of the vehicular traffic flow, and the changes in wind speed and direction, it is anticipated that there will be a complex distribution of air pollutants in both space and time. In this context, the utilization of mean values will be important to evaluate the influence of the urban landscape on the air pollution over the whole city,  but also in selected zones of particular interest.  Thus, we define the spatial mean values $\phi_k$ as the time functions:
\begin{equation}\label{eq:Mean_spacial_phi}
\phi_{k}(t) = \frac{1}{{\rm meas}(\Omega_{k})}\int_{\Omega_{k}} \phi(\mathbf{x},t) d\mathbf{x},
\end{equation}
and its corresponding time-averaged values,  defined as the constants:
\begin{equation}\label{eq:Mean_phi}
\Phi_k= \frac{1}{T}\int^T_0 \phi_{k}(t)dt,
\end{equation}
both computed in each selected urban zone $\Omega_k\subseteq\Omega$,  for $ k=1,\ldots,N_z$. 
%
\section{Numerical solution}\label{sec:NumSol}

In this section, we propose a numerical methodology for obtaining stable solutions to the aforementioned models. This methodology employs the finite element method of Lagrange $\rm{P_1}$ type, in conjunction with appropriate marching time schemes. 
As previously mentioned, all models are formulated within the same domain, i.e., the same mesh may be employed despite the differences between the spatial domains for the state variables. In this regard, let us denote the polygonal approximation of the domain as $\Omega_h$ and define a suitable triangulation $\mathcal{K}_h$ with $n_t$ triangles and $n_p$ nodes, where $n_t$ is the number of triangles and $n_p$ is the number of nodes. Consequently, we define over the polygonal approximation different functional spaces belonging to $P_1$ or $[P_1]^2$ on each triangle $K_i\in\mathcal{K}_h$ such that satisfies (depending on the model) all essential boundary conditions.  Concerning time discretization,  we take $N\in\mathbb{Z}^+$,  define the time step $\Delta t = T/N$,  and fix the following time instants $t^n = n\Delta t, \, n=0,\ldots,N$. 

Despite the unique characteristics of each numerical scheme, all of them are addressed in a consistent manner with respect to their full discrete problem. Specifically, the transition from the discrete variational form to a system of linear or nonlinear equations is enabled by the utilization of a suitable basis. In this study, the focus is on the final stage, and the references where further details can be found are provided for the interested reader.  
The initial example will be the traffic flow model and its desired speed, which is derived from the eikonal problem.

\subsection{Numerical solution of the traffic-flow model}

The solution of the nonconservative traffic-flow model (\ref{mod:traffic-flow})-(\ref{cond-ini:traffic-flow}), including the resolution of the eikonal problem (\ref{mod:eikonal}), has been extensively analyzed by the authors in their recent work \cite{garcia-chan2025}, where a full algorithm for its computation is presented. So, for the sake of conciseness, we will only recall here the main issues of its resolution, addressing interesting readers to \cite{garcia-chan2025} for further details.
In this case, we consider the functional spaces:
\begin{align*}
V_h = \{ v\in \mathcal{C}(\bar\Omega): v\vert_{K_i}\in P_1, \forall i = 1, \dots, n_t, \,\nabla v\cdot\mathbf{n} = 0 \mbox{ on } \Gamma_{in}\cup\Gamma_{out}\cup\Gamma_{w}\},\\
\mathbf{W}_h = \left\{\mathbf{w}\in [\mathcal{C}(\bar\Omega)]^2: \mathbf{w}\vert_{K_i}\in [P_1]^2, \forall i = 1, \dots, n_t, \,\nabla\mathbf{w}\mathbf{n}=\mathbf{0} \mbox{ on } \Gamma_{in}\cup\Gamma_{out},\right. \\ 
\left. \mbox{ and } \mathbf{w}\cdot\mathbf{n}=0 \mbox{ on } \Gamma_{w} \right\} ,
\end{align*}
for the finite element approximation of cars density $\rho$ and speed $\mathbf{u}$, respectively.  For the time semidiscretization, we use a high order strong stability preserving discretization \cite{gottlieb}, which is an explicit in time scheme of predictor-corrector type, useful to preventing spurious oscillations in nonlinear problems. In this way, after all the computations, we arrive to the discrete approximations $\{\{\rho^n_{h,i}\}^{n_p}_{i=1}\}^{N}_{n=0}$ and $\{\{\mathbf{u}^n_{h,i}\}^{n_p}_{i=1}\}^N_{n=0}$ for cars density and speed at each mesh node and each time instant. 


%
\subsection{Numerical solution of the air flow model}

As in above case, the authors have already introduced a full numerical algorithm for the resolution of problem (\ref{mod:airflow})-(\ref{bnd:airflow}) in previous papers \cite{garcia-chan2023,garcia-chan2024}, where the pressure stabilization method  \cite{chorin}  was applied. This method consists of splitting the momentum 
equation in two parts,  a first one for air velocity, and a second one for air pressure, using the constraint given by the continuity equation. In contrast with the methodology for the traffic flow model, the pressure stabilization method starts with an explicit time semidiscretization, followed by a standard finite element method, employing now the functional spaces: 
\begin{align*}
\hat V_h = \{ v\in \mathcal{C}(\bar\Omega) : v\vert_{K_i} \in P_1, \forall i = 1, \dots, n_t, \, \nabla v\cdot\mathbf{n} = 0 \mbox{ on } \Gamma_{in}\cup\Gamma_{w}\}, \\
\mathbf{\hat W}_h = \left\{\mathbf{w}\in [\mathcal{C}(\bar\Omega)]^2: \mathbf{w}\vert_{K_i}\in [P_1]^2, \forall i = 1, \dots, n_t, \nabla\mathbf{w}\mathbf{n}=\mathbf{0} \mbox{ on } \Gamma_{out},\right. \\ 
\left. \mbox{ and } \mathbf{w}\cdot\mathbf{n}=0 \mbox{ on } \Gamma_{w} \right\} ,
\end{align*}
for computing the discrete approximations $\{\{P^n_{h,i}\}^{n_p}_{i=1}\}^N_{n=0}$ and $\{\{\mathbf{v}^n_{h,i}\}^{n_p}_{i=1}\}^N_{n=0}$ for wind pressure and velocity, respectively.


\subsection{A stabilized finite element methodology to solve the air pollution model}

On the contrary to above cases, the numerical resolution of the air pollution model (\ref{mod:airpollution}) has not been previously addressed by our team. Thus, in this subsection we present the full details for obtaining the corresponding computational algorithm. 

The transport model (\ref{mod:airpollution}) is based in a classical advection-diffusion-reaction PDE with a dominant advective component. Utilizing the finite element method to solve this equation can result in undesirable spurious oscillations. Therefore, a stabilized methodology is required to avoid these unwanted effects, which have the propensity to overestimate the pollutant concentration and,  ultimately,  compromise the numerical solution.  In this paper, the Galerkin-Least-Square (GLS) methodology is used to stabilize the numerical solution. Therefore, we begin with the following semidiscrete in time problem. Let us evaluate the equation (\ref{eq:pollution}) at the time instant $t^{n+1}$,  and denote it as usual by $\phi^{n+1}$.  Now, let $v$ be a test function belonging to a suitable finite element space $\tilde V_h$. Multiplying equation (\ref{eq:pollution}) by $v$ and integrating on the domain $\Omega$ we have:
\begin{align}
\int\epsilon D_t{\phi}^{n+1}\,v + \int\mathbf{v}^{n+1}\cdot\nabla\phi^{n+1}\,v - \int\nabla\cdot\epsilon\nu\nabla\phi^{n+1}\,v + \int\epsilon\sigma\phi^{n+1}\,v =\nonumber\\ 
\int EC^{n+1}\,v ,
\end{align}
where $D_t$ denotes a discrete operator characterizing the time marching scheme, and the finite element space $\tilde V_h$ is defined as Lagrange polynomial functions associated to the triangular mesh $\mathcal{K}_h$:
\begin{equation}\label{space:pollution}
\tilde{V}_h = \{ v\in\mathcal{C}(\bar\Omega): v\vert_{K_i} \in P_1, \forall i = 1, \dots, n_t, \,\nabla v\cdot\mathbf{n} = 0 \mbox{ on } \Gamma_{out}\cup\Gamma_w\}.
\end{equation}
Then, after integrating by parts the diffusive term, we arrive to the variational form:
\begin{align}
\int\epsilon D_t{\phi}^{n+1}\,v + \int\mathbf{v}^{n+1}\cdot\nabla\phi^{n+1}\,v + \int\epsilon\nu\nabla\phi^{n+1}\cdot\nabla v  \nonumber\\ 
-\int_{\Gamma_{in}}\phi^{n+1} \mathbf{v}^{n+1}\cdot\mathbf{n} + \int\epsilon\sigma\phi^{n+1}\,v = \int EC^{n+1}v .
\end{align}
The GLS stabilizer consists of minimizing the residual $\mathcal{R}^{n+1}(\phi):= \mathcal{L}^{n+1}\phi + \epsilon\sigma\phi - EC^{n+1}$,  being the operator $\mathcal{L}^{n+1}$ defined from the transport equation 
(\ref{eq:pollution}) by $\mathcal{L}^{n+1}:= \mathbf{v}^{n+1}\cdot\nabla - \nabla\cdot\epsilon\nu\nabla$ (see \cite{larson2013,logg2012}). Thus,  we enunciate the following Least-Square scheme: find $\phi^{n+1}$ such that satisfies the following discrete problem:
\begin{align}
\int\epsilon D_t{\phi}^{n+1}v + \int\mathbf{v}^{n+1}\cdot\nabla\phi^{n+1} v + \int\epsilon\nu\nabla\phi^{n+1}\cdot\nabla v -\int_{\Gamma_{in}}\phi^{n+1} \mathbf{v}^{n+1}\cdot\mathbf{n}  \nonumber \\ 
+ \int\epsilon\sigma\phi^{n+1} v +\sum_{i=1}^{n_t}\int_{K_i} \mathcal{R}^{n+1}(\phi^{n+1})\Upsilon^{n+1}(v) = \int EC^{n+1} v,
\end{align}
for all $v\in \tilde{V}_h$, and where the operator $\Upsilon^{n+1}(v) = \zeta_1\mathbf{v}^{n+1}\cdot\nabla v - \zeta_2\nabla\cdot\epsilon\mu\nabla v$ is conformed by the streamline-diffusion and the artificial-diffusion stabilizer terms (the former adds diffusion,  but only in the wind direction, meanwhile the latter does it in all directions),  and modulated by the parameters $\zeta_1,\zeta_2>0$.  Then,  dropping the artificial-diffusion term by regularity requirements, and using a backward Euler marching time scheme we arrive to the following problem: find $\phi^{n+1}$ such that satisfies the following GLS problem:
\begin{align}\label{eq:GLS-full}
\int\epsilon\frac{\phi^{n+1}-\phi^n}{\Delta t}v + \int\mathbf{v}^{n+1}\cdot\nabla\phi^{n+1} v + \int\epsilon\nu\nabla\phi^{n+1}\cdot\nabla v -\int_{\Gamma_{in}}\phi^{n+1} \mathbf{v}^{n+1}\cdot\mathbf{n} \nonumber \\ 
+ \int\epsilon\sigma\phi^{n+1} v + \sum_{i=1}^{n_t} \int_{K_i} (\mathbf{v}^{n+1}\cdot\nabla\phi^{n+1}+\epsilon\sigma\phi^{n+1})(\mathbf{v}^{n+1}\cdot\nabla v) \nonumber \\
= \int EC^{n+1} v +\sum_{i=1}^{n_t} \int_{K_i} EC^{n+1} \mathbf{v}^{n+1}\cdot\nabla v.
\end{align} 
It is well known that introducing the nodal basis $\mathcal{B}_h = \{\tilde{v}_i \}^{n_p}_{i=1}$ of space $\tilde{V}_h$,  choosing $v=\tilde{v}_i$, and assuming the representation $\phi^{n+1}_i = \xi_i^{n+1}\tilde{v}_i$,  the problem (\ref{eq:GLS-full}) is equivalent to solve the following initial value problem:

Given the vector $\boldsymbol{\xi}^n= \{\xi^n_i\}^{n_p}_{i=1}$,  for each time instant $n=1,\ldots,N$, and given the initial state $\boldsymbol{\xi}^0=\{\phi^0_{i}\}^{n_p}_{i=1}$,  find the vector $\boldsymbol{\xi}^{n+1}$ (i.e.,  update the air pollution concentration) by solving the linear system: 
\begin{equation}\label{eq:linearsystem}
\frac{1}{\Delta t}\mathsf{M}_h (\boldsymbol{\xi}^{n+1} - \boldsymbol{\xi}^n) + [\mathsf{C}_h^{n+1} + \mathsf{R}_h - \mathsf{Bnd}^{n+1}_{in} + \hat{\mathsf{M}}_h + \mathsf{SL}^{n+1}_h]\boldsymbol{\xi}^{n+1} = \mathsf{f}_h^{n+1} + \mathsf{sl}^{n+1}_h 
\end{equation}
where the matrices and vectors could depend on time via the wind velocity, and are obtained after the standard assembling process of the finite element method.  For the sake of conciseness,  we omit here the list of integrals that yield each matrix and vector of equation (\ref{eq:linearsystem}) (e.g. $(\mathsf{M}_h)_{i,j} = \int \epsilon \tilde{v}_j\tilde{v}_i$ and $(\mathsf{f}^{n+1}_h)_i = \int EC^{n+1}\tilde{v}_i$), but it is possible to identify them easily by comparing equations (\ref{eq:linearsystem}) and (\ref{eq:GLS-full}).

Above numerical methodology has shown enough good properties to avoid spurious oscillations, as can be seen in below section.   

\begin{rmk}
It is worthwhile remarking here that above GLS scheme (\ref{eq:GLS-full}) is a simplified version,  due to the omission of stabilizer terms related to the time marching scheme $D_t\phi^{n+1}$ on the residual $\mathcal{R}$ (see \cite{logg2012} for further details).  Nevertheless, the application of this stabilizer streamline-diffusion method, as successfully employed by authors in \cite{garcia-chan2023,garcia-chan2024},  has already been demonstrated effective in addressing a similar heat transport problem on an urban-porous medium. 
\end{rmk}
%
\section{Numerical experiments} \label{sec:NumExp}
In this section, we present numerical results for the different models, including a thorough explanation of the mesh processing and the utilization of a common mesh for the traffic flow,  the air flow,  and the pollution transport models.  Thus, the space-time distribution of the cars density and their local velocity are presented,  and according to the selection of carbon dioxide $\rm{CO_2}$ as the pollution indicator, we proceed also to show the space distribution of instantaneous emissions by vehicle and the emissions concentration.  Finally,  we generate a reference wind field to simulate air pollutant transport.  Results are presented and discussed here as an element of the overall phenomena of urban air pollution.  The domain chosen for our computational examples corresponds to the Guadalajara Metropolitan Area,  one of the largest metropolitan zones in Mexico.

\subsection{Definition of the study zone,  and characterization of obstacles and urban/rural zones for the models}
\begin{figure}
\centering
\subcaptionbox{\label{fig1:GoogleEarth}}{\includegraphics[width=.45\linewidth,height=.4\textheight,keepaspectratio]{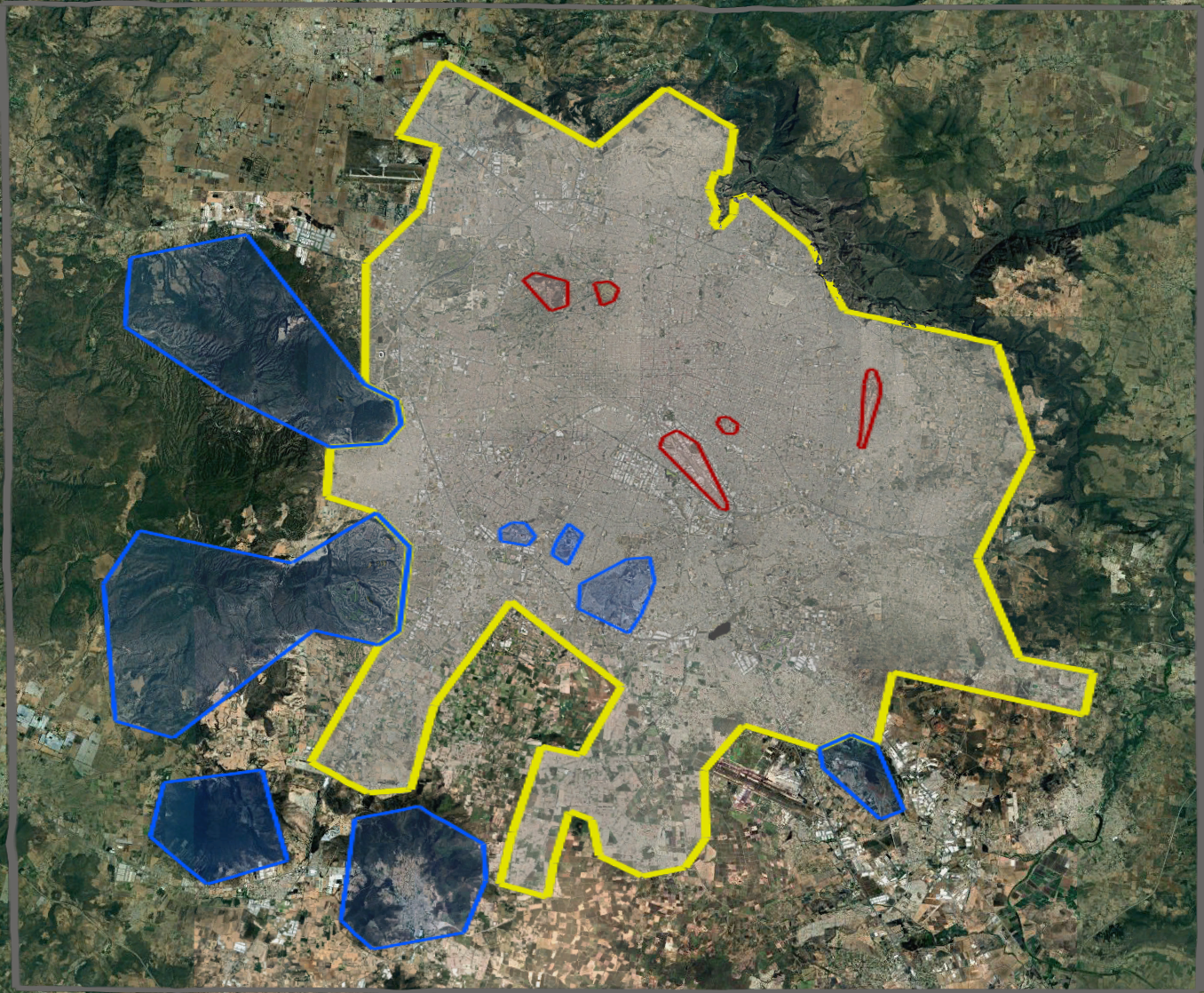}}
\subcaptionbox{\label{fig1:GDL-mesh}}{\includegraphics[width=.45\linewidth,height=.6\textheight,keepaspectratio]{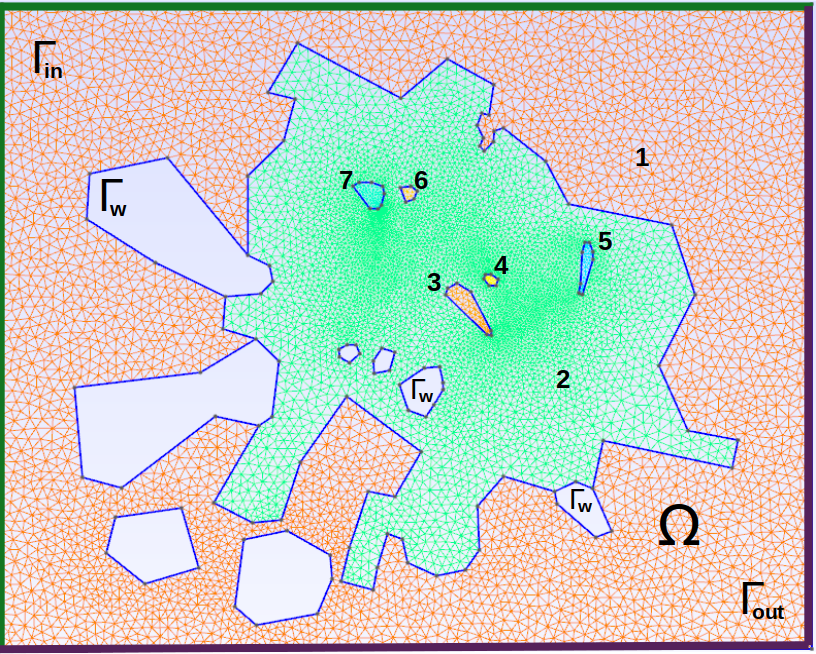}}
\caption{(a) The Metropolitan Area of Guadalajara and its rural surroundings,  taken from Google-Earth (2024),  and (b) the triangular mesh generated with the software Gmsh,  showing the numbers for the different zones and the labels of the boundary segments.}\label{fig1}
\end{figure}

Our study is based on a domain $\Omega$ which encompasses the city and its rural surroundings within a rectangular area.  With respect to the elements inside the rectangle, we have the city delimited by a number of straight lines conforming the city layout.  A few natural obstacles (for instance, mountains) are also idealized with regular shapes.  Some of these obstacles are situated within the urban zone, while others are located beyond its layout. 
In Figure \ref{fig1:GoogleEarth},  these natural  obstacles correspond to the areas in blue,  the urban zone corresponds -once excluded all above natural obstacles- to the area inside the yellow polygonal line,  and the rural zone corresponds to the area outside it.
Additionally, our domain includes several selected zones within the city that exhibit low urbanization and are restricted to traffic flow.  Despite the numerous zones that could impede traffic flow, a selection was made to represent all the possible types.  The five selected areas are shown in Figure \ref{fig1:GoogleEarth} delimited by red lines.  They are numbered in Figure \ref{fig1:GDL-mesh} as follows: (3) industrial zone,  (4) university campus, (5) Solidaridad park,  (6) Golf Club,  and (7) Colomos park.  In Figure \ref{fig1:GDL-mesh},  numbers (1) and (2) stand for the rural zone and the urban one, respectively.

In order to ensure the applicability of the same mesh or parts of the same mesh to all the mathematical models in this paper, the domain must be meshed with care.  Consequently, each geographical element previously referenced (rural zone, urban zone, obstacles, and selected areas) is delineated herein as a polygon, with each polygon nested within the other ones.  Specifically, the rural zone surrounds the urban zone and a number of natural obstacles, while the selected areas and the rest of the natural obstacles are located within the urban zone.  The Gmsh software has been selected for the mesh generation process \cite{gmsh2009}.  Gmsh is capable of meshing each polygon in a nested manner, which is initiated with the rural zone (meshed in orange), followed by the urban zone (meshed in green),  and concluded with the selected areas (meshed in diverse colors: red, yellow, blue...).  However, depending on the model under study,  there is the option to either omit these selected areas for the traffic flow model or to mesh them for the wind and the pollution transport models. The decision of meshing or not meshing the selected areas is made without altering the numbering in the remaining domain (urban and rural zones).  Consequently,  for the traffic flow model, the mesh has fewer elements (9232 nodes and 17893 triangles) compared to the air flow and the air pollution models, which have 9363 nodes and 18258 triangles (this latter mesh can be seen in Figure \ref{fig1:GDL-mesh}).

\begin{figure}
\centering
\subcaptionbox{\label{fig2:dense-city}}{\includegraphics[width=.45\linewidth,height=.6\textheight,keepaspectratio]{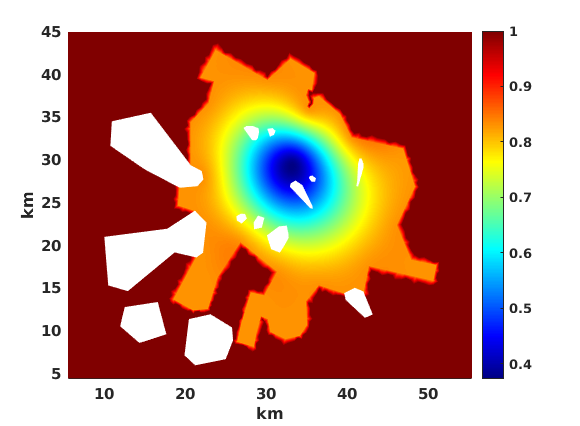}}
\subcaptionbox{\label{fig2:disperse-city}}{\includegraphics[width=.45\linewidth,height=.6\textheight,keepaspectratio]{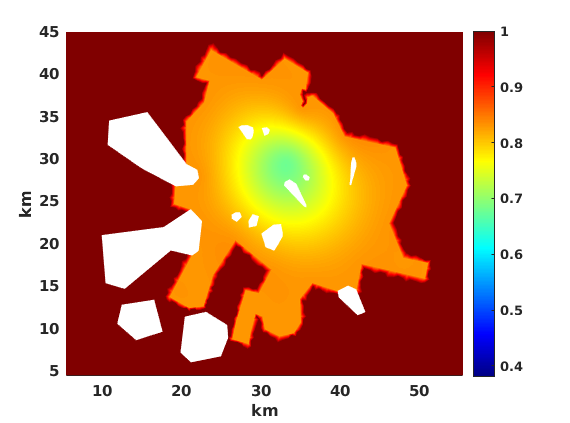}}
\caption{Porosity distribution in the domain $\Omega$ for the two city configurations: (a) The dense buildings city,  with a porosity value of 0.38 at the center and a value of 0.82 on the city limits.  (b) The disperse city, showing a 0.68 porosity value at its center and similar values to the dense case on suburbia.}\label{fig2:city-porosity}
\end{figure}

In order to address the concept of urban porosity, this paper considers two distinct city configurations based on a concentric city model: the dense city and the disperse city.  The former corresponds to the case of a city with high density of buildings (and, consequently, low porosity),  and the latter to a low building density with higher porosity. To establish the porosity values for both cities, a porosity value is designated for the city's downtown area, and another value is assigned to the suburbs,  following the extension of porosity over the entire city via interpolation. The values for urban porosity were set between 0.38 and 0.82, as described in \cite{hu2012}. Consequently, the dense city is characterized by a porosity value of 0.38 in its central area, while the disperse city exhibits there a porosity value of 0.68. In both scenarios, the city layout demonstrates a porosity value of 0.82. Conversely, the rural zone is characterized as empty (without buildings), with a porosity value of 1.0 (see Figure \ref{fig2:city-porosity}).
 
Concerning the boundary conditions, we divided the boundary domain into three segments: the inlet part $\Gamma_{in}$, the outlet part $\Gamma_{out}$, and the natural obstacles part $\Gamma_{w}$. So, depending on the model, we use these boundary segments to impose different conditions on the traffic flow, the wind flow, and the air pollution transport (see Figure \ref{fig1:GDL-mesh}).

\subsection{The traffic flow simulations and their associated pollutant emissions}

\begin{figure}
\centering
\subcaptionbox{\label{fig3:init-cond}}{\includegraphics[width=.45\linewidth,height=.6\textheight,keepaspectratio]{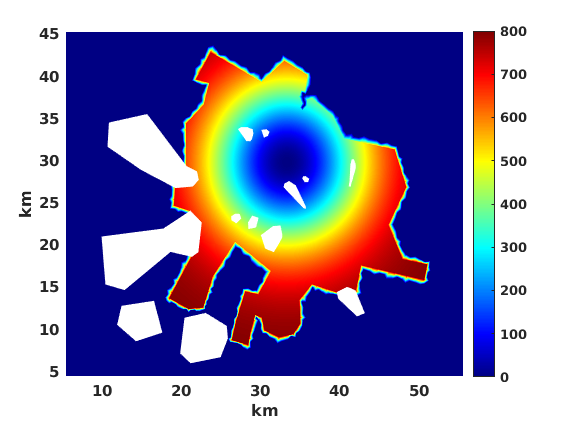}}
\subcaptionbox{\label{fig3:absorption}}{\includegraphics[width=.45\linewidth,height=.6\textheight,keepaspectratio]{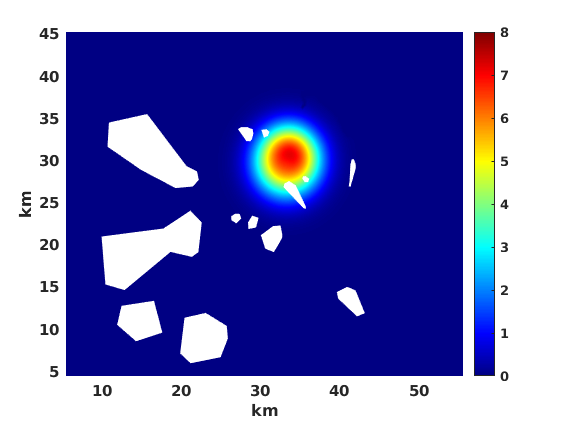}}
\subcaptionbox{\label{fig3:eikonal}}{\includegraphics[width=.45\linewidth,height=.6\textheight,keepaspectratio]{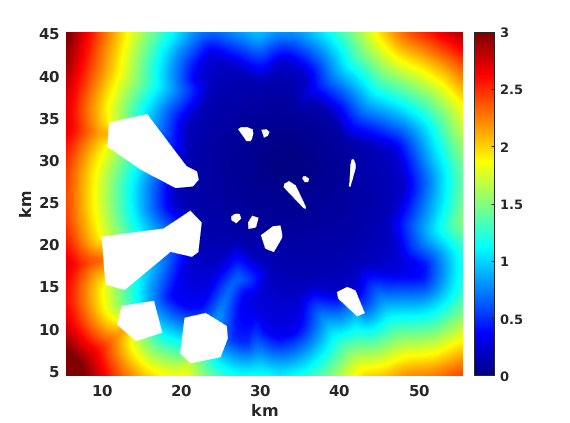}}
\subcaptionbox{\label{fig3:desired-direction}}{\includegraphics[width=.45\linewidth,height=.6\textheight,keepaspectratio]{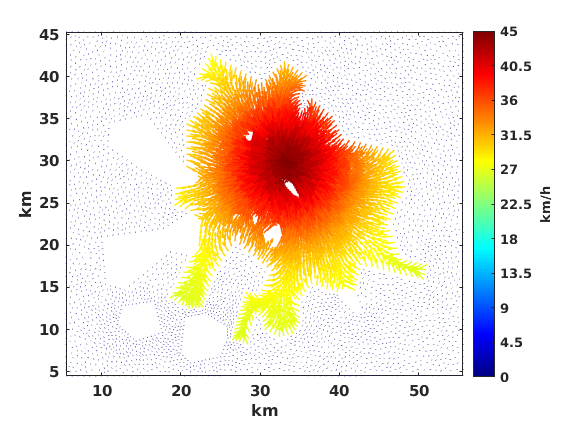}}
\caption{(a) Initial condition for the vehicular density.  (b) Absorption rate of the solid phase (off-street parking) for the dense city.  (c) Distribution of the eikonal potential $\varphi$ on the dense city.  (d) Desired traffic speed at initial time $t=0$.}\label{fig3:init-cond-desired-velocity}
\end{figure}

In this section, we present the traffic flow dynamics describing the vehicular density and its local speed as solutions of the traffic flow model (\ref{mod:traffic-flow}). We simulate the traffic flow over 2 hours with a simple pattern: at the initial time cars are dispersed concentrically by a Gaussian distribution with center at the downtown and suitable variance (i.e.,  traffic demand is $q=0$), so high cars density values are present on suburbia contrasting with an almost empty city center.  When time begins to march, drivers flow to their destinations in the city downtown where use off-street parking inside buildings to place their cars,  leaving free space on the streets for other cars. This pattern is performed until all the cars leave the streets which is plausible due to an unlimited parking capacity.  So, in Figures \ref{fig3:init-cond} and \ref{fig3:absorption} we show the initial state of cars density and the street-off parking rate,  noting in the latter how cars can leave the domain only going out through the downtown area. The desired velocity $\mathbf{u}_d$ is the driver's reference to reach its destination,  being generated by the eikonal solution $\varphi$,  which is displayed in Figure \ref{fig3:eikonal} for the first instant simulation.  Notice that eikonal values inside the city are lower and in particular $\varphi=0$ at the center.  The desired speed and direction are shown in Figure \ref{fig3:desired-direction}, mainly pointing to the city center where the maximum velocity $U_{max} = 45\,\rm{km/h}$ is reached.  In suburbia, an interesting phenomenon happens: where the urban part is a narrow zone surrounded by rural areas, the desired direction targets a path that allows cars to leave this area first,  to retake later the path to the center (that is, it looks like cars driving to take a main road avoiding the rural area).  Meanwhile, in the rural zone, the desired velocity is almost zero,  targeting outside the city,  but with a null influence on vehicles behavior.  

\begin{figure}
\centering
\subcaptionbox{\label{fig4:density-densecity-0pt5h}}{\includegraphics[width=.45\linewidth,height=.6\textheight,keepaspectratio]{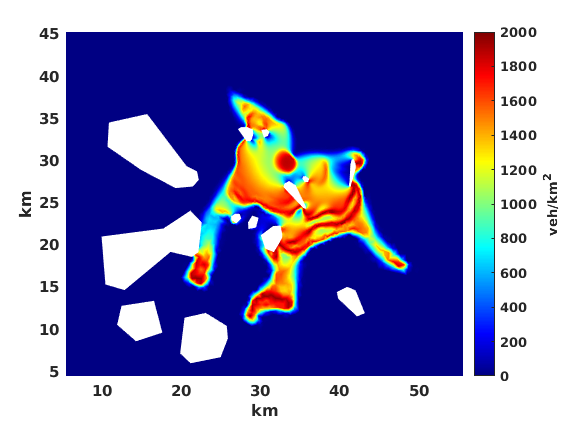}}
\subcaptionbox{\label{fig4:density-dispersecity-0pt5h}}{\includegraphics[width=.45\linewidth,height=.6\textheight,keepaspectratio]{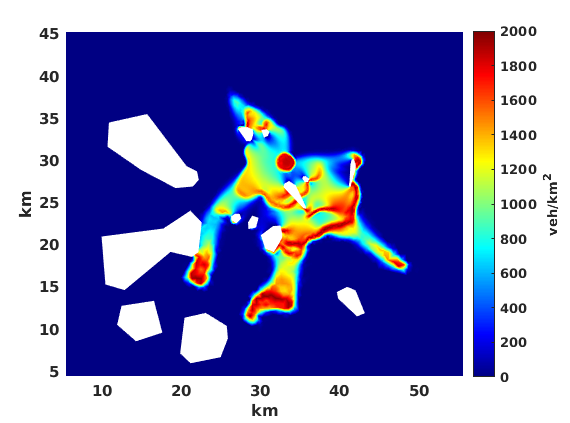}}
\subcaptionbox{\label{fig4:density-densecity-1pt0h}}{\includegraphics[width=.45\linewidth,height=.6\textheight,keepaspectratio]{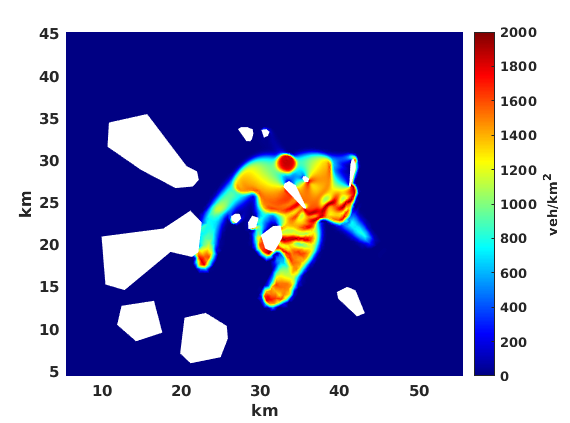}}
\subcaptionbox{\label{fig4:density-dispersecity-1pt0h}}{\includegraphics[width=.45\linewidth,height=.6\textheight,keepaspectratio]{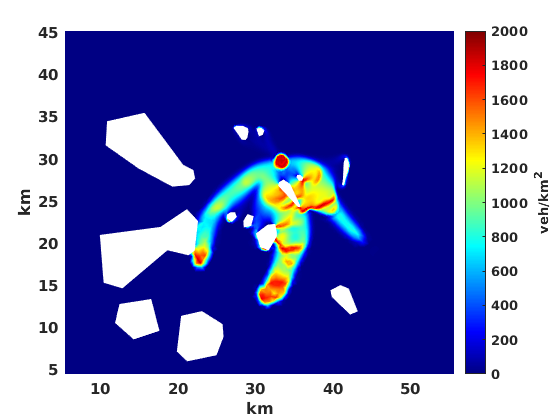}}
\subcaptionbox{\label{fig4:density-densecity-1pt5h}}{\includegraphics[width=.45\linewidth,height=.6\textheight,keepaspectratio]{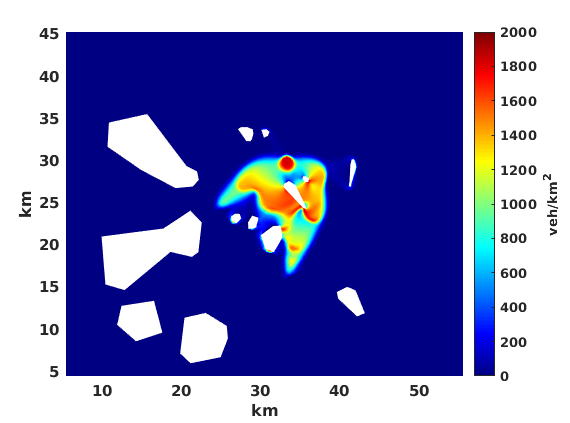}}
\subcaptionbox{\label{fig4:density-dispersecity-1pt5h}}{\includegraphics[width=.45\linewidth,height=.6\textheight,keepaspectratio]{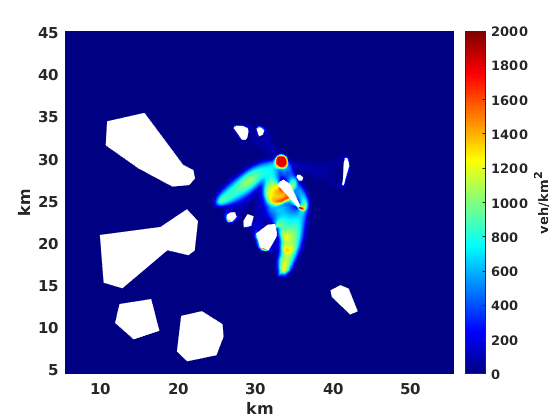}}
\caption{Vehicular density $\rho$ $[\rm{veh/km^2}]$ for the dense city (left column) and the disperse city (right column) at times  0.5, 1.0, 1.5 h (from top to bottom).}\label{fig4:density}
\end{figure}

The desired velocity is defined as an ideal scenario in an urban environment devoid of architectural obstructions. It is anticipated that drivers will deviate from this ideal velocity by adjusting their direction and speed according to the available space for flow and the prevailing vehicle density.  As illustrated in Figure \ref{fig4:density}, we show frames depicting the density evolution at times 0.5, 1.0, and 1.5 hours (from top to bottom) both for a dense city (left column) and a disperse city (right column).  It is evident that traffic jams are present in both cases and at all times due to an insufficient street-off parking rate,  to the narrow passes, and to obstacles.  However, in the dense city, these jams extend over a larger area and persist over time, requiring more time for the city configuration to evacuate all vehicles. Additionally, rarefaction waves, where density achieves maximum values but with low velocities, are present due to drivers in front traveling rapidly to occupy free space produced by the off-street parking, and drivers at the wave back accelerating until reach the wave peak.  These wave dynamics are attributable to the incorporation of a traffic flow model that takes into account the correction proposed by \cite{Aw2000},  in order to address the limitations of macroscopic models identified by \cite{daganzo1995}.  In this way, drivers are primarily motivated by information available to them at the front of their vehicles.    

\begin{figure}
\centering
\subcaptionbox{\label{fig5:speed-densecity-0pt5h}}{\includegraphics[width=.45\linewidth,height=.5\textheight,keepaspectratio]{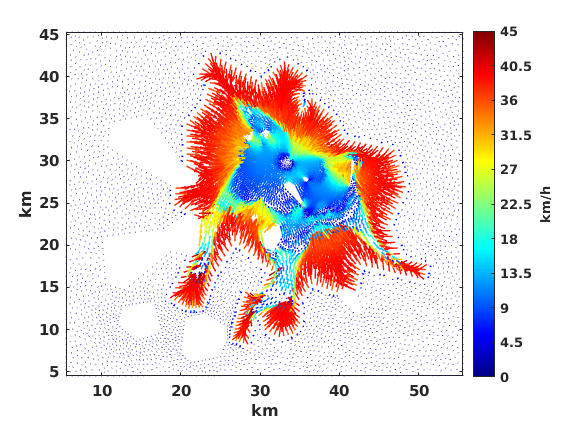}}
\subcaptionbox{\label{fig5:speed-dispersecity-0pt5h}}{\includegraphics[width=.45\linewidth,height=.5\textheight,keepaspectratio]{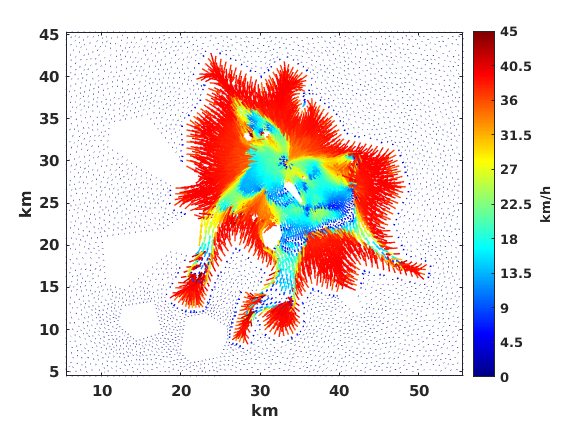}}
\subcaptionbox{\label{fig5:speed-densecity-1pt0h}}{\includegraphics[width=.45\linewidth,height=.5\textheight,keepaspectratio]{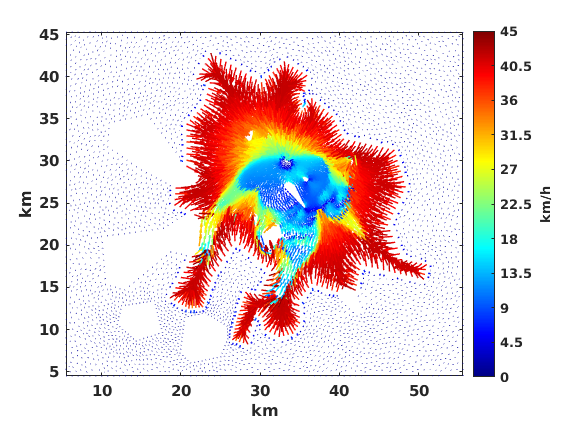}}
\subcaptionbox{\label{fig5:speed-dispersecity-1pt0h}}{\includegraphics[width=.45\linewidth,height=.5\textheight,keepaspectratio]{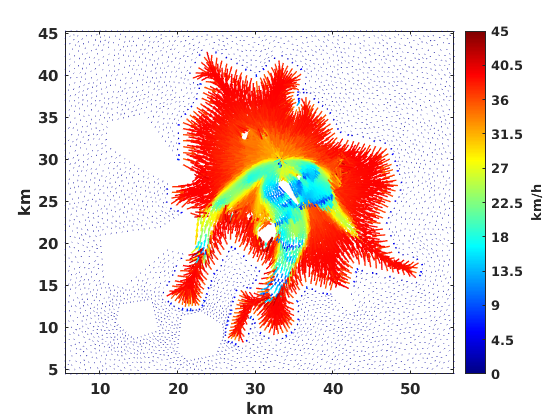}}
\subcaptionbox{\label{fig5:speed-densecity-1pt5h}}{\includegraphics[width=.45\linewidth,height=.5\textheight,keepaspectratio]{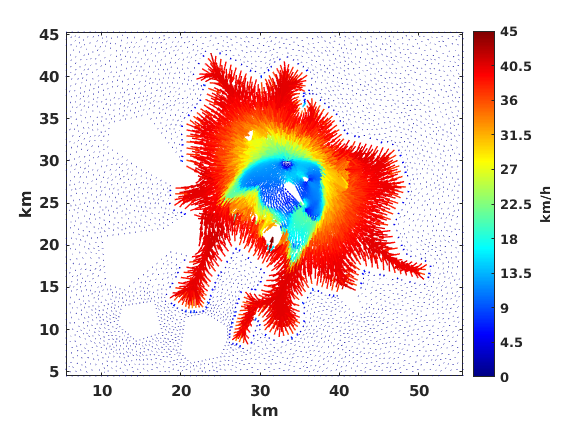}}
\subcaptionbox{\label{fig5:speed-dispersecity-1pt5h}}{\includegraphics[width=.45\linewidth,height=.5\textheight,keepaspectratio]{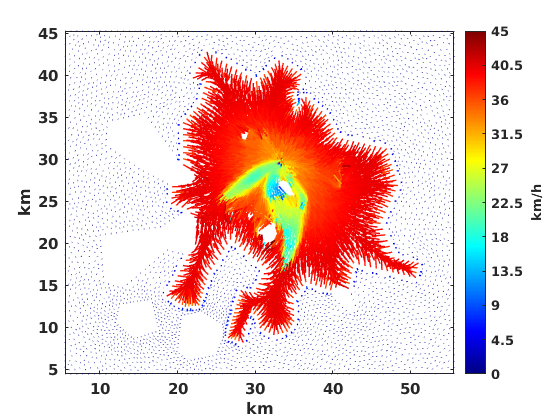}}
\caption{Traffic local velocity $\mathbf{u}$ $[\rm{km/h}]$ for the dense city (left column) and the disperse city (right column) at times  0.5, 1.0, 1.5 h (from top to bottom).}\label{fig5:speed}
\end{figure}

\begin{figure}
\centering
\subcaptionbox{\label{fig6:acceleration-densecity}}{\includegraphics[width=.45\linewidth,height=.5\textheight,keepaspectratio]{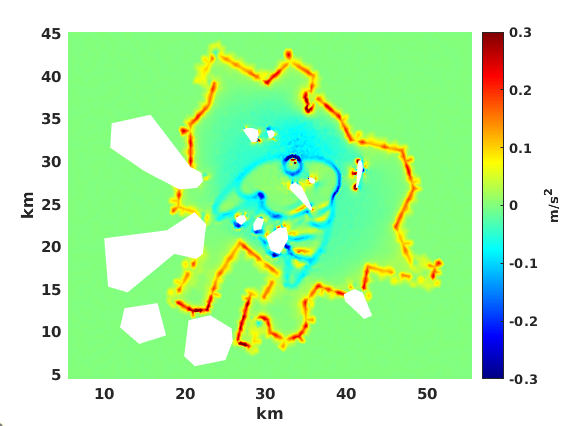}}
\subcaptionbox{\label{fig6:acceleration-dispersecity}}{\includegraphics[width=.45\linewidth,height=.5\textheight,keepaspectratio]{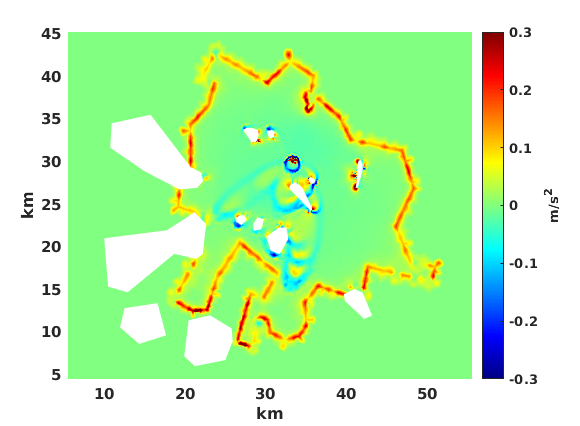}}
\caption{Traffic flow acceleration $a$ $[\rm{km/h^2}]$ computed from equation (\ref{eq:scalar-acceleration}),  at time 1.5 h,  for the dense city (a) and the disperse city (b).}\label{fig6:acceleration}
\end{figure}

Concerning traffic flow,  local speed depends on the desired traffic velocity and the porous media dynamics (mainly, the Darcy Law and its corrections).  Therefore,  differences between both city configurations are expected.  Indeed, the disperse city presents higher velocity values,  allowing drivers to reach the center in less time,  facilitating an ordered parking, and showing a less congested downtown.  By contrast,  the dense city complicates traffic flow due to the high density of buildings generating less velocity values.  So,  drivers need more time to reach the center,  generating more car congestions.  As the density,  the local speed field presents also waves,  and faces the obstacles with a tangential direction to overrun them.  All of these above issues can be observed in Figure \ref{fig5:speed}.  With respect to car acceleration, both cities show very similar patterns.  In the downtown area, characterized by lower porosity values,  a deceleration can be noticed.  In the rarefaction waves,  a typical acceleration-deceleration dynamics is found.  Finally, on the obstacle walls, the acceleration behavior is noteworthy: when vehicles encounter an obstacle, they decelerate until the obstacle is behind them; at that point, they accelerate to reach their previous velocity (see Figure \ref{fig6:acceleration} for both types of city).  

\begin{figure}
\centering
\subcaptionbox{\label{fig8:ECO2-densecity}}{\includegraphics[width=.45\linewidth,height=.5\textheight,keepaspectratio]{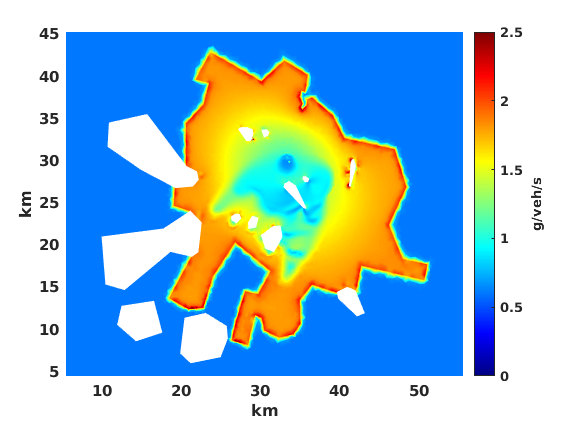}}
\subcaptionbox{\label{fig8:ECO2-dispersecity}}{\includegraphics[width=.45\linewidth,height=.5\textheight,keepaspectratio]{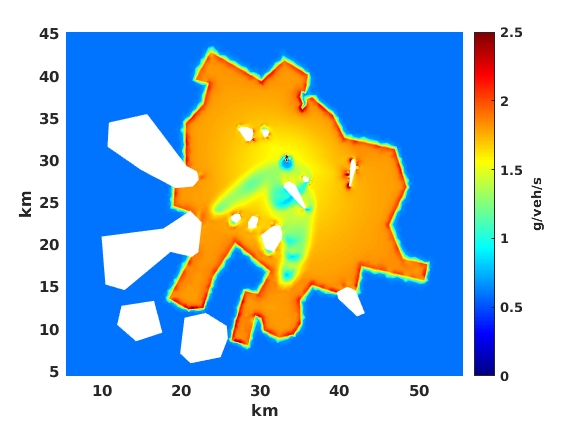}}
\caption{Instantaneous emissions $E_{\rm{CO_{2}}}$ [g/veh/s] for the dense city (a) and disperse city (b),  computed from equation (\ref{eq:instant-emissions}) at time 1.5 h.} 
\end{figure}

\begin{figure}
\centering
\subcaptionbox{\label{fig9:densecity-0pt5h}}{\includegraphics[width=.45\linewidth,height=.5\textheight,keepaspectratio]{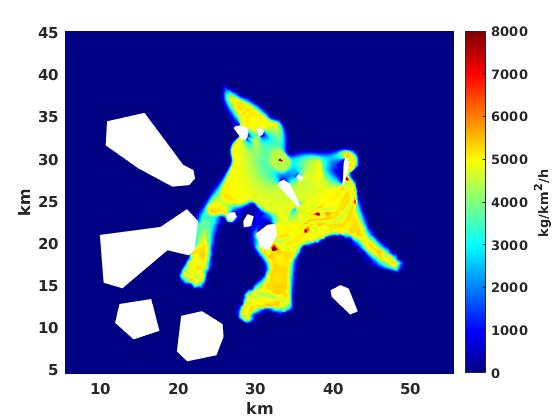}}
\subcaptionbox{\label{fig9:dispersecity-0pt5h}}{\includegraphics[width=.45\linewidth,height=.5\textheight,keepaspectratio]{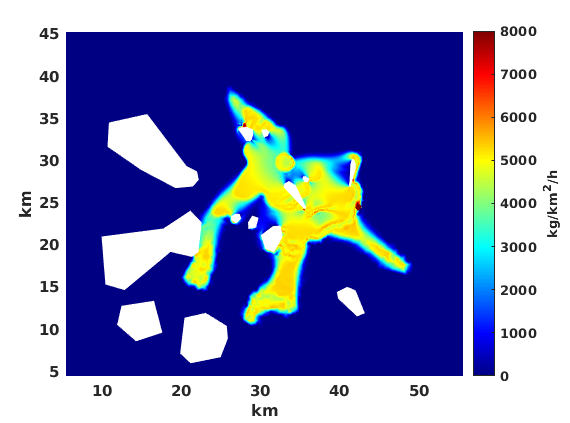}}
\subcaptionbox{\label{fig9:densecity-1pt0h}}{\includegraphics[width=.45\linewidth,height=.5\textheight,keepaspectratio]{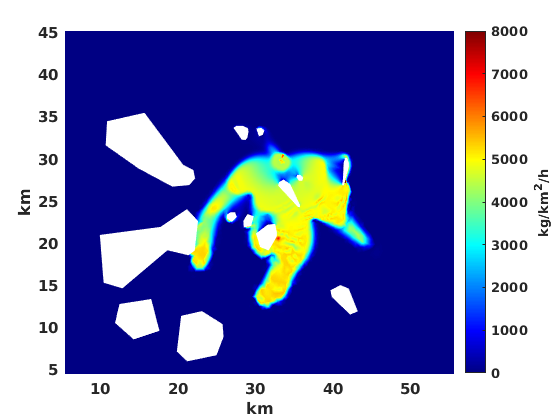}}
\subcaptionbox{\label{fig9:dispersecity-1pt0h}}{\includegraphics[width=.45\linewidth,height=.5\textheight,keepaspectratio]{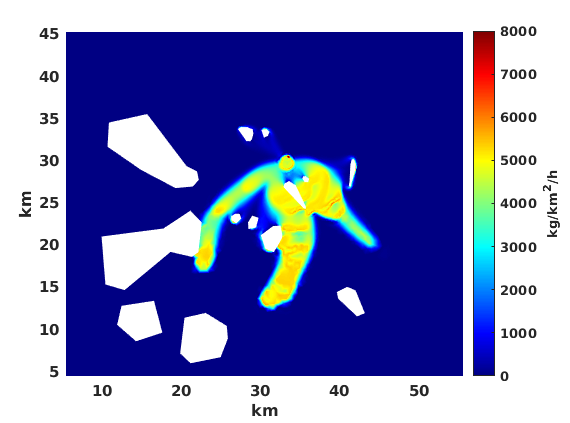}}
\subcaptionbox{\label{fig9:densecity1pt5h}}{\includegraphics[width=.45\linewidth,height=.5\textheight,keepaspectratio]{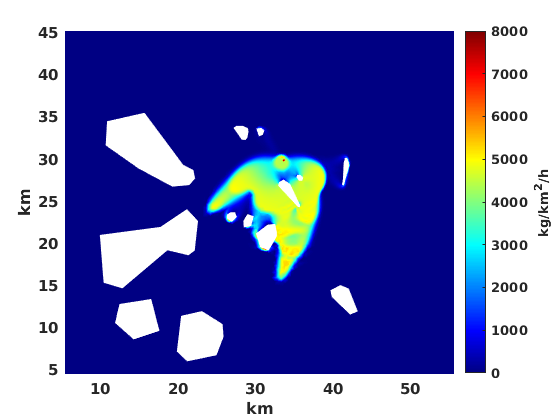}}
\subcaptionbox{\label{fig9:dispersecity-1pt5h}}{\includegraphics[width=.45\linewidth,height=.5\textheight,keepaspectratio]{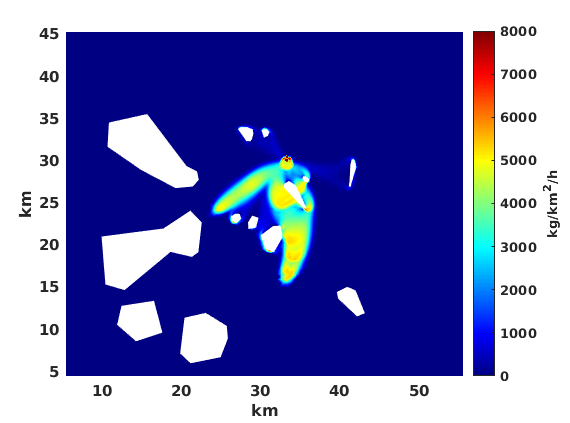}}
\caption{Spatial distributions of the emissions concentration $EC$ of $\rm{CO_{2}}$ [kg/km$^2$/h] computed using equation (\ref{eq:concent-emissions}) for the dense city (left column) and the disperse city (right column) at times 0.5, 1.0, 1.5 h (from top to bottom).}\label{fig9:emissions-concentration}
\end{figure}

The emissions are computed using the microscopic emissions model (\ref{eq:instant-emissions}) along with an equation for distributed emissions concentration (\ref{eq:concent-emissions}). The inputs for this model are two scalar functions (one for velocity $U$ and one for acceleration $a$) computed from equations (\ref{eq:scalar-velocity}) and (\ref{eq:scalar-acceleration}),  depending on previously computed solutions $\mathbf{u}$ and $\mathbf{a}$.  In this paper, we focus only on $\rm{CO_{2}}$,  despite the microscopic emissions model could manage others (PM,  $\rm{NO_{x}}$,  VOC and so on).  Based on \cite{panis},  $\rm{CO_2}$ microscopic emissions are directly related to velocity, so lower velocity agrees with lower $\rm{CO_2}$ emissions but, at the same time, acceleration facilitates $\rm{CO_2}$ emissions.  So, in Figures \ref{fig8:ECO2-densecity} and \ref{fig8:ECO2-dispersecity} we show the $\rm{CO_2}$ distribution at 1.5 hours.  In both cases,  when comparing these Figures to Figures \ref{fig5:speed-densecity-1pt5h} and \ref{fig5:speed-dispersecity-1pt5h} for velocity and Figures \ref{fig6:acceleration-densecity} and \ref{fig6:acceleration-dispersecity} for acceleration,  we can observe two important issues.  First,  where a large number of vehicles is present (in the downtown area) they cause low velocities,  which draw $\rm{CO_2}$ emissions to low levels.  Then,  in contrast,  high instantaneous emissions close to almost empty zones are given by acceleration (generating at some areas a combination of high/low emissions depending on acceleration/deceleration).  Also,  it is worthwhile remarking here that instantaneous emissions correspond to each car (are measured in [g/veh/s]),  but multiplying them by density we get the emissions concentration (in [kg/km$^2$/h]),  as shown in Figure \ref{fig9:emissions-concentration}.

Emissions concentration is directly comparable to vehicles density,  as can be confirmed by checking Figures \ref{fig4:density} and \ref{fig9:emissions-concentration}. Highly polluted areas due to traffic congestion, including the top of rarefaction waves,  the circumvented obstacles,  the selected areas and the downtown zone,  are identified as particularly problematic.  The emissions persist within the city until all vehicles are removed by off-street parking.  The configuration of the city is a key factor in understanding the emissions patterns.  In a dense city,  cars congestions are more prevalent,  affecting larger urban areas,  but emissions are lower due to the reduced velocity of vehicles.  Conversely, in a disperse city,  congestions are limited to smaller areas, but emissions are higher due to the increased velocity of vehicles.  This leads to the conclusion that both dense and disperse cities can exhibit areas with low or high levels of $\rm{CO_2}$ emissions.  However, a clear advantage of the disperse city is that cars can be evacuated more quickly,  which implies a emissions concentration lower than in the dense city.

\subsection{A reference wind field}

\begin{figure}
\centering
\subcaptionbox{\label{fig7:wind-norm-densecity}}{\includegraphics[width=.45\linewidth,height=.6\textheight,keepaspectratio]{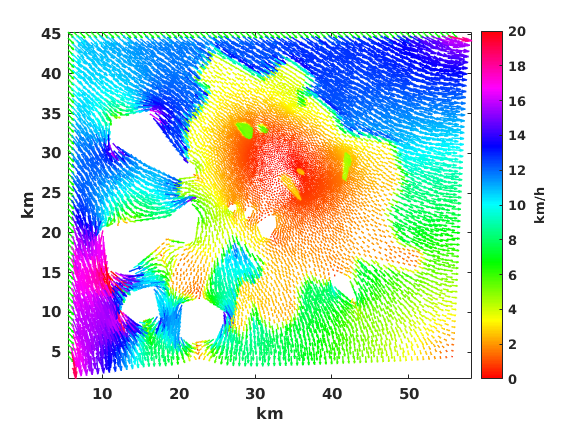}}
\subcaptionbox{\label{fig7:wind-norm-dispersecity}}{\includegraphics[width=.45\linewidth,height=.6\textheight,keepaspectratio]{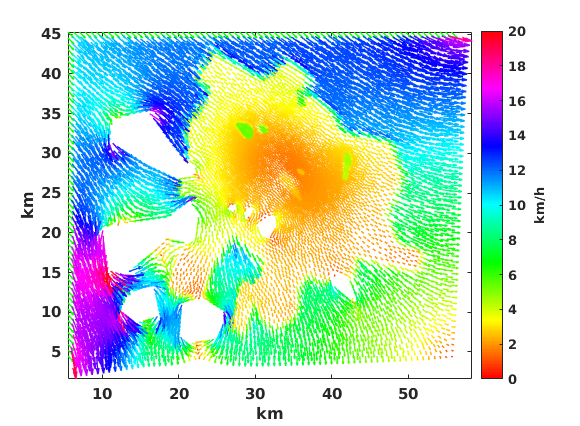}}
\caption{Wind velocity $\mathbf{v}$ $[\rm{km/h}]$ norm and direction for the dense city (a) and the disperse city (b).}\label{fig7:reference-wind}
\end{figure}

We propose an idealized representation of the complex wind field within the atmospheric boundary layer, employing a constant reference wind field. This reference wind field is derived by imposing a constant inlet wind with a northwestern directional component, specified as $\mathbf{u}_{in}$ = (5,-5) $\rm{km/h}$ (or,  equivalently,  (1.38,-1.38) $\rm{m/s}$) on the inlet section of the boundary $\Gamma_{in}$.  After a sufficiently large time, the numerical solution of the air flow model reaches a stationary state, which is used here as our reference wind field (see Figure \ref{fig7:reference-wind}). The reference wind field in the rural area exhibits velocities within the range of $0$ to $20$ km/h. The abnormally high velocities in some small regions are,  in this case,  a direct consequence of the narrow pass due to our particular domain geometry (see the bottom left corner of Figures \ref{fig7:wind-norm-densecity} or \ref{fig7:wind-norm-dispersecity}).  Additionally,  the slip condition (\ref{bnd:airflow-slip-wall}) imposed on walls facilitates the air flow to circumvent obstacles (see Figure \ref{fig7:reference-wind}).  However, in the urban area,  both city configurations show notable differences,  as can be noted in Figure \ref{fig7:wind-norm-densecity}.  In the dense city,  the wind velocity norm is close to zero due to the narrow space between buildings (low porosity) in the center,  but reaches more than 4 km/h in suburbia due to wider streets and less area covered with buildings (high porosity).  In contrast, the disperse city exhibits greater variability in wind velocities, with values ranging from 2 to 4 km/h across the entire urban zone (see Figure \ref{fig7:wind-norm-dispersecity}).  This fact suggests that air pollutants could be more readily evacuated in a disperse city,  as opposed to a dense city where stagnant air pollutants are likely to accumulate due to insufficient space for air circulation.  In the selected zones, that exhibit high porosity values with respect to their urban surroundings, wind velocities approximate those observed in rural or suburban areas, reaching 8 km/h in the parks and 4 km/h in the industrial zone and the university campus (another notable feature of our modeling approach).  This phenomenon can be attributed to the contrasting porosity values of the model, which delineate a constricted pathway within the urban landscape.  Consequently, the wind velocity within these selected zones is elevated to ensure the satisfaction of the incompressibility restriction imposed by the continuity equation (\ref{eq:trafficflow-continuity}).

\subsection{Air pollution due to traffic flow}
\begin{figure}
\centering
\subcaptionbox{\label{fig10:densecity-0pt5h}}{\includegraphics[width=.45\linewidth,height=.5\textheight,keepaspectratio]{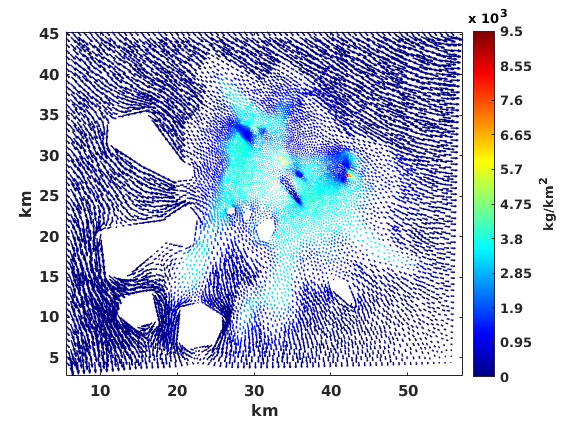}}
\subcaptionbox{\label{fig10:disperecity-0pt5h}}{\includegraphics[width=.45\linewidth,height=.5\textheight,keepaspectratio]{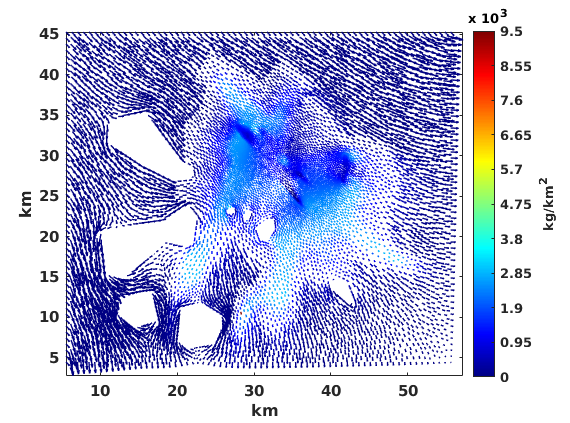}}
\subcaptionbox{\label{fig10:densecity-1pt0h}}{\includegraphics[width=.45\linewidth,height=.5\textheight,keepaspectratio]{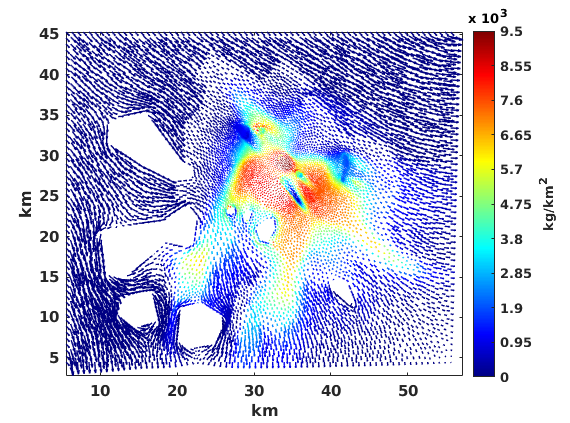}}
\subcaptionbox{\label{fig10:disperecity-1pt0h}}{\includegraphics[width=.45\linewidth,height=.5\textheight,keepaspectratio]{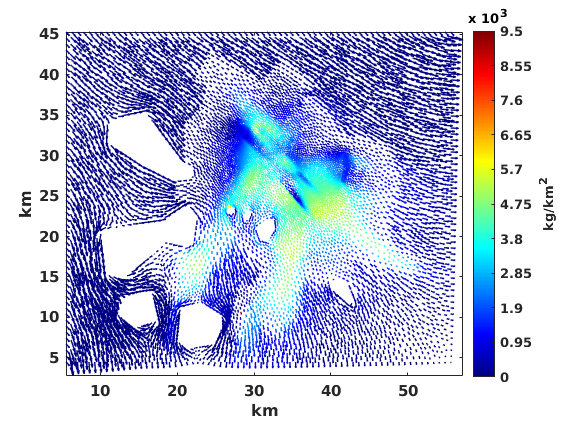}}
\subcaptionbox{\label{fig10:densecity-1pt5h}}{\includegraphics[width=.45\linewidth,height=.5\textheight,keepaspectratio]{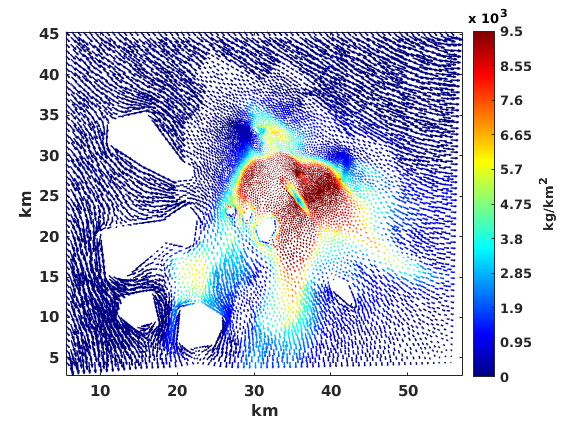}}
\subcaptionbox{\label{fig10:disperecity-1pt5h}}{\includegraphics[width=.45\linewidth,height=.5\textheight,keepaspectratio]{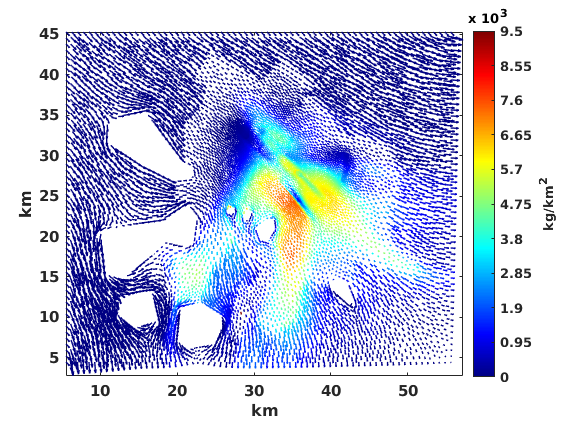}}
\caption{Air pollution concentration $\phi$ $[\rm{kg/km^2}]$ due to emissions of $\rm{CO_{2}}$ for the dense city (left column) and the disperse city (right column) at times 0.5, 1.0, 1.5 h (from top to bottom).}\label{fig10:air-pollution-CO2}
\end{figure}

In this section, we will examine the impact of traffic on air pollutant evacuation efficiency for a city,  considering its built environment.  The air flow results suggest that a disperse city may have superior conditions for polluted air evacuation due to its higher air speeds compared to a dense city.  As illustrated in Figure \ref{fig10:air-pollution-CO2}, the left column depicts the air pollution concentration (in [$\rm{kg/km^2}$]) for a dense city,  while the right column showcases the corresponding values for a disperse city.  The temporal progression of air pollution is also delineated from top to bottom, with values at 0.5, 1.0, and 1.5 hours.  It is noteworthy that the disparities between the two types of cities intensify over time.  At 0.5 hours, both cities exhibit comparable air pollution distribution patterns, with elevated concentrations in downtown and urban-confined areas, while rural zones in the right, left and bottom directions of the urban zone experience comparatively lower levels.  However, in the dense city, pollution levels are slightly higher.  These dynamics intensify at 1.0 h, with a more pronounced deterioration in the downtown area of the dense city, while the disperse city experiences more extensive pollution plumes but with lower concentrations.  Finally,  at 1.5 hours, a significant proportion of vehicles in both cities have already reached their designated parking areas.  This results in increased wind activity, which now plays the main role in air circulation and pollutant dissipation.  At this point, the advantages of the dispersed city become apparent, characterized by better air circulation, a more expansive but limited pollution plume, and a downtown area in the northwestern section that is comparatively uncontaminated.  In contrast, in the dense city, the pollution plumes are extensive, spanning almost the entire downtown area,  and exhibiting much higher levels of pollutants.  This observation suggests that the dispersed city may experience faster pollutant dissipation.

\begin{figure}
\centering
\subcaptionbox{Urban area\label{fig11:mean-urban}}{\includegraphics[width=.45\linewidth,height=.6\textheight,keepaspectratio]{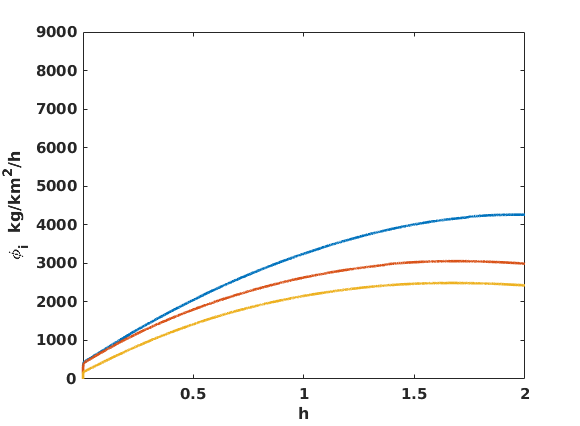}}
\subcaptionbox{Dense city\label{fig11:mean-densecity}}{\includegraphics[width=.45\linewidth,height=.6\textheight,keepaspectratio]{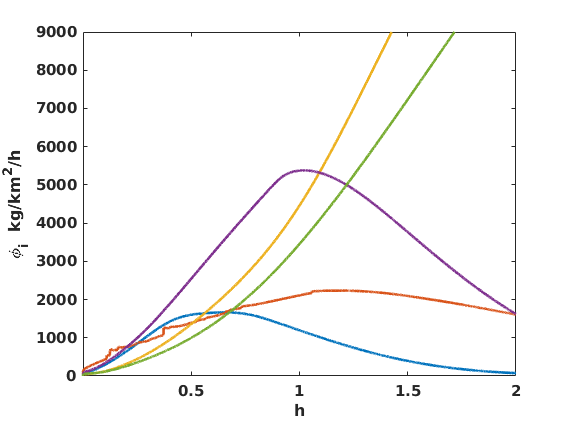}}
\subcaptionbox{Disperse city \label{fig11:mean-dispersecity}}{\includegraphics[width=.45\linewidth,height=.6\textheight,keepaspectratio]{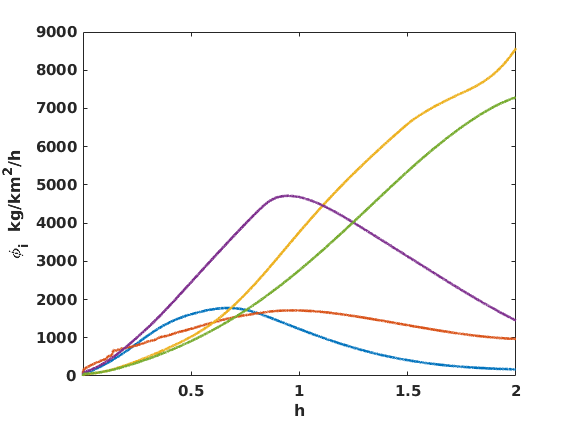}}
\subcaptionbox{Disperse city with $U_{max} = 30$ km/h\label{fig11:mean-dispersecity-limitspeed}}{\includegraphics[width=.45\linewidth,height=.6\textheight,keepaspectratio]{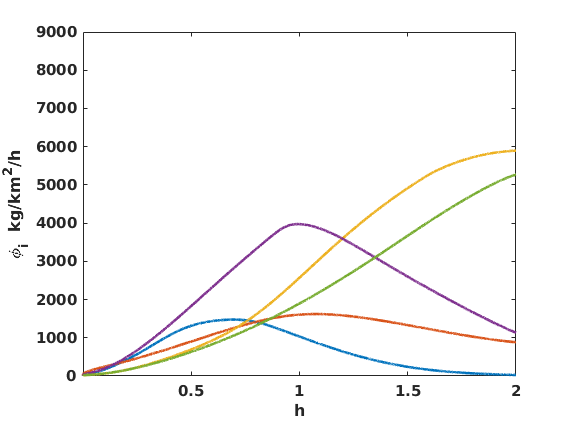}}
\caption{Time evolution of the mean air pollution concentration $\phi_{\textcolor{blue}{k}}$ [kg/km$^{2}$] for the whole urban zone (a) and for the selected zones (b)-(d).  In (a) we show the mean values for the dense city (blue), the disperse city (red),  and the disperse city with lower speed limit (yellow).  In (b)-(d), the five curves correspond to Colomos park (blue), Solidaridad park (red), university campus (yellow), Golf Club (purple),  and industrial zone (green).}\label{fig11:mean-air-pollution}
\end{figure}

Notably,  Figure \ref{fig10:air-pollution-CO2} also shows variances in air pollution concentrations across both space and time dimensions, thereby enabling the identification of areas with acceptable or unacceptable air quality levels, needed for a better management of the air pollution.  To further analyze these dynamics,  Figure \ref{fig11:mean-air-pollution} presents the evolution of the mean concentration (computed from expression (\ref{eq:Mean_spacial_phi})) over time, both in the city as a whole and in specific urban zones.  Thus, in Figure \ref{fig11:mean-urban} we show three curves corresponding to three cases: a dense city, a disperse city, and a disperse city where maximum velocity $U_{max}$ of vehicles is limited to only 30 km/h, instead of previous 45 km/h.  It is clear that the dense city represents the worst case, but it is also shown that limiting the speed to 30 km/h is the best case.  Concerning selected zones, it is interesting to observe that porosity does not influence in the same manner all selected zones.  Indeed, in all the configurations of the city, the areas far from the city center, such as Colomos and Solidaridad parks, and the Golf Club show a bell-shaped evolution of their mean concentration, increasing until reaching a maximum but showing after this a significant reduction.  On the contrary,  in the selected zones close to the center, such as the university campus and the industrial zone, the mean concentration increases throughout the whole simulation (see Figures \ref{fig11:mean-densecity}-\ref{fig11:mean-dispersecity-limitspeed}).

Thus,  the Golf Club, that presents its worst situation in the dense city,  improves in the case of the dispersed city,  but its best situation is reached by imposing the speed limit.  The Colomos park, after reaching its maximum mean concentration of $\rm{CO_2}$,  drops it to minimum levels at the end of the simulation in the three city configurations.  However, Solidaridad park did not demonstrate a notably reduction of its levels after reaching its maximum value,  as was observed in Colomos park, being this tendency even slightly worse in the disperse city.  Regarding the university campus and the industrial zone, their tendency is always increasing, changing only the intensity of this increase,  being worse in the dense city but improving notably with the imposition of speed limits (see Figures \ref{fig11:mean-densecity}-\ref{fig11:mean-dispersecity-limitspeed}).

In order to provide further support for the findings presented in Figure \ref{fig11:mean-air-pollution}, we show in Table \ref{tab1:mean-values} the average $\rm{CO_2}$ concentrations, calculated using equation (18),  for the whole city but also for three the selected zones.  From Table 1 we can see that the worst case of the average concentration is represented by the dense city, both in the whole city and in the selected zones.  On the other hand, the disperse city presents a reduction in all zones with the exception of Colomos park,  with a slightly higher mean than the one obtained in the dense city.  This may be attributable to the high traffic speed reached in the zone close to the park,  which requires greater acceleration after circumventing the park in the case of the dispersed city.  Finally,  it is confirmed that limiting speed reduces the average $\rm{CO_2}$ concentration in all cases.

Above results indicate that limiting vehicles speed reduces $\rm{CO_2}$ concentration in all city configurations and all selected zones.  
This may suggest that a mitigation strategy consisting of the choice of a speed limit depending on the sensibility of the zone (i.e.,  a different speed limit for the different areas of the city) could be an effective and easy to implement solution to reduce air pollution.


\begin{table}[b]
\centering
\begin{tabular}{lcccc}
Type of city             & Urban zone & Colomos & University & Solidaridad \\
\hline
Dense city       & 2.9522$\times 10^{-3}$             &  0.8360$\times 10^{-3}$              & 5.6995$\times 10^{-3}$                  & 1.6697$\times 10^{-3}$   \\
Disperse city    & 2.3299$\times 10^{-3}$             &  0.8730$\times 10^{-3}$              & 3.8441$\times 10^{-3}$                  & 1.2506$\times 10^{-3}$   \\
Limited speed  & 1.8717$\times 10^{-3}$             &  0.6628$\times 10^{-3}$              & 2.7833$\times 10^{-3}$                  & 1.1008$\times 10^{-3}$   \\
\hline
\end{tabular}
\caption{Mean values of air $\rm{CO_2}$ concentration $\Phi_k$ $[\rm{kg/km^2}]$ on the whole urban zone and some of the selected zones (Colomos and Solidaridad parks,  and the university campus) for the dense city, the disperse city and the limited speed case.}\label{tab1:mean-values}
\end{table}

%
\section{Conclusions}\label{sec:Conclusions}

In this paper, a multi-model study of urban air pollution due to a nonconservative traffic flow in a large porous city was made from a numerical viewpoint.  Several PDEs models were formulated and combined to model the key elements involved in this important environmental issue.  Thus,  four different models for traffic flow, air flow, microscopic emissions, and air pollutant transport were studied and adapted to our porous city scenario.  Numerical schemes to solve these models were also established briefly: the traffic flow model is solved using a standard finite element methodology,  the air flow model uses the pressure stabilization method,  and the air pollutant transport model was solved with a Least-Square stabilization method.  Concerning time marching schemes, an explicit type was used for traffic flow and air flow, but an implicit one was applied in the case of the transport model.

Numerical simulations suggest the direct influence of the urban landscape in all the elements of the phenomenon.  Indeed,  porosity influences the traffic flow speed and acceleration, the air flow magnitude,  the emissions by car and the total emissions, and, consequently,  this also influences the air pollution transport.  In particular, the dense cities,  where buildings occupy more urban surface,  are related to a lower traffic flow speed with changes in acceleration located near obstacles and in the downtown.  In consequence, the dense city presents low instantaneous emissions in most of the city,  showing notorious changes in the city layout and near obstacles.  The spatial distribution of emissions is more persistent in time covering more urban areas -due to cars needing more time to make off-street parking in the city downtown.  So,  due to low airflow speed,  dense cities present less capacity to evacuate air pollution.

In contrast, the disperse city,  where cars have more space to flow, is related to higher traffic flow velocities (consequently,  producing higher emissions of pollutants). Nevertheless, with a higher velocity,  cars can evacuate the city in less time,  decongesting the downtown.  Concerning to evacuation of pollutants, as expected, the disperse city facilitates it with a higher air flow speed.

In summary,  we highlight two final important results: the proposal of speed limit as a potential air pollution mitigation tool in our porous city,  and the capability of our two-dimensional model to reproduce the main aspects of the phenomena involved in the process,  including rarefaction waves typical in one-dimensional models.

%
\section*{CRediT authorship contribution statement}
\textbf{N\'estor Garc\'ia-Chan:} Conceptualization, Methodology, Software, Investigation, Writing – original draft \textbf{Lino J. Alvarez-V\'azquez:} Conceptualization, Methodology, Investigation, Writing – review and editing \textbf{A. Mart\'inez:} Conceptualization, Methodology, Investigation, Writing – review and editing \textbf{Miguel E. V\'azquez-M\'endez:} Conceptualization, Methodology, Investigation, Writing – review and editing
\section*{Declaration of competing interest}
The authors declare that they have no known competing financial interests or personal relationships that could have appeared to influence the work reported in this paper.
\section*{Acknowledgements}
This work was supported by the SNII-CONAHCYT and PRODEP (Mexico) [grant numbers SNII-52768, 103/16/8066]; Ministerio de Ciencia e Innovación (Spain) and NextGenerationEU (European Union) [grant number TED2021-129324B-I00].



\bibliographystyle{elsarticle-num-names} 







\end{document}